\documentclass[preprint,12pt,number,sort&compress]{elsarticle}
\journal{Advances in Mathematics}
\usepackage[T1]{fontenc}
\usepackage{lmodern}
\usepackage{amsmath,amssymb,amsthm,mathtools}
\usepackage{booktabs,longtable,tabularx,array}
\usepackage{adjustbox}
\usepackage{enumitem}
\usepackage{microtype}
\usepackage{xcolor}
\usepackage{hyperref}
\hypersetup{colorlinks=true,linkcolor=blue!45!black,citecolor=blue!45!black,urlcolor=blue!45!black}
\newtheorem{theorem}{Theorem}[section]
\newtheorem{proposition}[theorem]{Proposition}
\newtheorem{lemma}[theorem]{Lemma}

\theoremstyle{definition}

\theoremstyle{remark}

\newcommand{\PP}{\mathbb P}
\newcommand{\CC}{\mathbb C}
\newcommand{\C}{\mathbb C}

\newcommand{\ii}{\mathrm i}

\newcommand{\ord}{\operatorname{ord}}
\newcommand{\Res}{\operatorname{Res}}

\newcommand{\calK}{\mathcal K}
\newcommand{\CCfam}[1]{\mathcal C_{#1}}

\renewcommand{\arraystretch}{1.16}
\makeatletter
\def\ps@pprintTitle{%
  \def\@oddhead{}%
  \def\@evenhead{}%
  \def\@oddfoot{\hfil\thepage\hfil}%
  \let\@evenfoot\@oddfoot
}
\makeatother
\begin{document}
\begin{frontmatter}
\title{Complete Center Classification of a Quintic $\mathbb Z_2$-Equivariant System via Rational Descent and Ramification}
\author[aff1]{Feng Li}
\ead{lifeng@lyu.edu.cn}
\author[aff1]{Meilan Cai}
\ead{caimeilan@lyu.edu.cn}
\author[aff1]{Chaochao Sun\corref{cor1}}
\ead{sunuso@163.com}
\cortext[cor1]{Corresponding author.}
\address[aff1]{School of Mathematics and Statistics, Linyi University, Linyi, China}
\begin{abstract}
We develop a new structural method to completely classify the elementary centers of a ten-parameter quintic $\mathbb Z_2$-equivariant system. The method overcomes the computational bottleneck of direct focus-quantity elimination by combining a Li\'enard reduction, a necessary-and-sufficient analytic involution criterion, and rational descent with residue and ramification analysis. We obtain forty-two explicit original-parameter conditions covering the entire center locus, including all degenerate and boundary cases. Sufficiency is established constructively through integrating factors, first integrals and explicit reversible realizations.
\end{abstract}
\begin{keyword}
Center--focus problem \sep quintic $\mathbb Z_2$-equivariant systems \sep Li\'enard equation \sep analytic involution \sep L\"uroth theorem \sep rational differential \sep ramification
\MSC[2020] Primary 34C07 \sep Secondary 34C05 \sep 34C20
\end{keyword}
\end{frontmatter}

\section{Introduction and main results}\label{sec:intro}
The center--focus problem is a central local question in the qualitative theory of planar polynomial differential equations. It is closely related to, but distinct from, the second part of Hilbert's sixteenth problem, which asks for bounds on the number and relative positions of limit cycles of polynomial vector fields \cite{Hilbert1902,Poincare1885,Li2003}. For a singular point with purely imaginary, nonzero eigenvalues, a \emph{center} is surrounded by a punctured neighborhood of periodic orbits, whereas a \emph{focus} has spiraling nearby trajectories. The Poincar\'e return map encodes this distinction: the center condition amounts to the vanishing of its displacement germ. Its coefficients, commonly called focus quantities or Lyapunov constants, have been studied since Poincar\'e and Lyapunov \cite{Poincare1885,Lyapunov1892}. Their simultaneous vanishing is an infinite analytic condition, while finite algebraic descriptions become possible in polynomial parameter families. Related Abel-equation formulations of the center--focus problem provide another useful perspective \cite{Gavrilov2020}.

The quadratic case provides the classical model for this interaction between analysis and algebra. Early work of Dulac, Kapteyn and Frommer developed the local center theory \cite{Dulac1908,Kapteyn1912,Frommer1934}. Bautin's theorem then showed that the local center and small-amplitude limit-cycle bifurcation questions for quadratic polynomial fields are controlled by a finite sequence of focal quantities \cite{Bautin1952}. The algebraic viewpoint has since been developed through Bautin ideals, reduced indices, perturbation methods and the geometry of center varieties \cite{Christopher2005,ChristopherLi2007,HauserRislerTeissier1999,Zoladek1994,DumortierRoussarieRousseau1994}. Importantly, the local cyclicity statement for quadratic systems does not determine the global Hilbert number for quadratic vector fields.

For cubic and higher-degree systems, the proliferation of parameters makes direct focal-value elimination considerably harder. Representative constructions of small-amplitude limit cycles and studies of perturbations can be found in \cite{JamesLloyd1991,Zoladek1995,YuHan2005,LiLiuYang2009,HanYu2012}. In parallel, the complex singular-point-quantity method and its recurrence formulas have provided a systematic basis for symbolic calculations \cite{LiuLi1990,Liu2001,LiuLiHuang2008,LiuHuang2005}. These techniques are powerful for detecting obstructions to a center and for investigating local bifurcations, but a long sequence of necessary conditions does not by itself constitute an exhaustive center classification. For an exact classification, the proof must also describe all exceptional parameter strata and verify sufficiency.

Symmetry offers additional structure. If a finite group $G$ acts linearly on the plane, a vector field $X$ is $G$-equivariant when $X(gz)=gX(z)$ for every $g\in G$. The \emph{action} of the group matters: in this paper $\mathbb Z_2$ acts by the central inversion $z\mapsto-z$, so that $X(-x,-y)=-X(x,y)$. This must not be confused with equivariance under a reflection of a coordinate axis. Central inversion interchanges the two opposite nonzero equilibria, which motivates the simultaneous, or bi-center, question. The geometric and bifurcation structures of equivariant polynomial fields have been investigated in \cite{LiLiu2010,LiChanChung2002,ChanChungLi2001,LiChanChung2002Z6,LiuLi2009Z5,LiuLi2009Z6}. In the cubic $\mathbb Z_2$ setting, the works \cite{YuHan2004Z2,YuHan2005Z2,LiuLi2011,LiuLi2014} address focal-value bifurcations and the bi-center problem. Related questions about isochronicity, linearization and nilpotent singularities require additional hypotheses and are treated, for other families, in \cite{DuLiu2017,LiLiuLiuYu2020,LiLiuLiuYu2018}; nonsmooth equivariant switching systems form yet another distinct class \cite{GuoYuChen2019}.

For quintic systems, studies of symmetric limit-cycle configurations and bifurcation sets \cite{LiChanChung2002,LiZhangMiao2012,KongXiong2013} demonstrate the richness of the degree-five problem. Other results concern selected equivariant quintic families, simultaneous centers, critical-period bifurcations, nilpotent singularities or integrability under specific invariant-line assumptions \cite{Romanovski2017,Gine2018,Fernandes2017,WuLi2023,LiWuChenYu2025,LiZhangLi2026}. Such findings are relevant background, but they do not automatically give a necessary-and-sufficient classification for the full normalized family considered below. Likewise, a result about centers should not be interpreted as an isochronicity theorem: the analytic coordinate changes and positive time changes used here need not preserve the period function.

We consider the real normalized family consisting of a linear field and homogeneous quintic terms,
\begin{equation}\label{eq:family}
\begin{aligned}
\dot x={}&-\frac\delta4x-(a_1+1)y+\frac\delta4x^5
+a_1x^4y+a_2x^3y^2+a_3x^2y^3+a_4xy^4+a_5y^5,\\
\dot y={}&-\frac14x+(\delta-a_6)y+\frac14x^5
+a_6x^4y+a_7x^3y^2+a_8x^2y^3+a_9xy^4+a_{10}y^5.
\end{aligned}
\end{equation}
The equilibria $(1,0)$ and $(-1,0)$ are interchanged by central inversion. At $(1,0)$, the Jacobian is $\bigl(\begin{smallmatrix}\delta&-1\\1&\delta\end{smallmatrix}\bigr)$, so $\delta=0$ is necessary for a center. The classification concerns this \emph{ten-parameter linear-plus-homogeneous-quintic model}, rather than all quintic $\mathbb Z_2$-equivariant fields with arbitrary cubic terms, and it addresses elementary centers rather than nilpotent or isochronous centers.

For the full ten-parameter family \eqref{eq:family}, a direct focus-quantity calculation encounters a substantial computational bottleneck. Shifting the equilibrium $(1,0)$ to the origin produces terms of several degrees, and successive focal quantities can become large multivariate polynomials. Their simultaneous elimination must retain every vanishing leading coefficient and every degenerate parameter branch; a long generic calculation alone therefore does not establish complete coverage.

These difficulties motivate the structural approach developed here. A local Li\'enard reduction first turns the center problem into a necessary-and-sufficient common-involution criterion. Differential identities then eliminate the nonrational integrating factor and yield rational invariants. Their common rational function field, together with residue and ramification constraints, reduces the classification to finitely many cases. We reconstruct the original coefficients on every retained parameter stratum and establish sufficiency independently by integrating factors, first integrals or changes to reversible systems. The contribution lies in combining these ingredients into an exhaustive classification, including the degenerate and boundary cases.

The resulting classification is stated directly in the original parameters.

\begin{theorem}[Complete center classification]\label{thm:main}
For the normalized real family \eqref{eq:family}, the equilibrium $(1,0)$ is an elementary center if and only if $\delta=0$ and the ten-tuple $(a_1,\ldots,a_{10})$ satisfies at least one of the conditions $\CCfam1,\ldots,\CCfam{42}$ given explicitly in \ref{app:atlas}. In that case $(-1,0)$ is also an elementary center by central inversion. The listed conditions cover every real exceptional and boundary parameter stratum; the list is a \emph{cover}, not a claim that its members are disjoint or irreducible components.
\end{theorem}

The forty-two conditions are grouped according to the rational-descent reconstruction: four elementary/degenerate conditions $\CCfam1$--$\CCfam4$, nine degree-two conditions $\CCfam5$--$\CCfam{13}$, sixteen degree-three conditions $\CCfam{14}$--$\CCfam{29}$, nine degree-four conditions $\CCfam{30}$--$\CCfam{38}$, and four degree-five conditions $\CCfam{39}$--$\CCfam{42}$. The degree here refers to the common rational generator, not to the quintic degree of the original vector field. The full original-parameter formulas are recorded exactly once in \ref{app:atlas}; the proof is independent of how the coefficient atlas is typeset.

The structural notation and branch labels in the classification atlas are explained in detail in the necessity analysis in Section~\ref{sec:classification}, using the definitions developed in Sections~\ref{sec:lienard}--\ref{sec:residue}.

The proof of Theorem~\ref{thm:main} proceeds in two parts. Sections~\ref{sec:lienard}--\ref{sec:residue} develop the analytic criterion and the rational tools needed for necessity. Section~\ref{sec:classification} applies them to obtain the complete parameter cover and summarizes the origins of the conditions in Table~\ref{tab:atlas}. Section~\ref{sec:sufficiency} proves sufficiency in the same numerical order as the conditions in
\ref{app:atlas}.

\section{Li\'enard reduction and the analytic center criterion}\label{sec:lienard}
Impose $\delta=0$ and work locally near $(1,0)$, where $x>0$. Introduce
\begin{equation}\label{eq:AB}
\begin{gathered}
s=\frac yx,\qquad r=x^4,\qquad c=a_1+1,\qquad g=a_6,\\
p(s)=\sum_{j=1}^5 a_j s^j,\qquad
q(s)=\frac14+\sum_{j=1}^5 a_{j+5}s^j,\\
A(s)=q(s)-sp(s),\qquad B(s)=-\frac14-gs+cs^2.
\end{gathered}
\end{equation}
In particular,
\begin{equation*}
A(s)=\frac14+gs+(a_7-a_1)s^2+(a_8-a_2)s^3
+(a_9-a_3)s^4+(a_{10}-a_4)s^5-a_5s^6.
\end{equation*}
Direct calculation yields
\begin{equation}\label{eq:sr}
\dot s=B+rA,\qquad \dot r=4r(-cs+rp).
\end{equation}
Set $v=\dot s$. Since $A(0)=1/4$, the inverse $r=(v-B)/A$ is analytic and $v$ is an actual local coordinate. Differentiating \eqref{eq:sr} gives
\begin{equation*}
\dot s=v,\qquad \dot v=h(s)v^2+f(s)v+\gamma(s),
\end{equation*}
where
\begin{equation}\label{eq:hfg}
\begin{split}
h&=\frac{A'+4p}{A},\\
f&=B'-\frac{BA'}A-4cs-\frac{8Bp}{A},\\
\gamma&=4Bcs+\frac{4B^2p}{A}.
\end{split}
\end{equation}
The notation $\gamma$ is used for the scalar term to avoid confusion with the coefficient $g=a_6$. The Jacobian of $(x,y)\mapsto(s,v)$ at $(1,0)$ is nonsingular. No parameter division has yet been used.

Define the normalized analytic functions
\begin{equation*}
\begin{aligned}
K(s)&=\exp\left(-\int_0^s h(t)\,dt\right),\\
F(s)&=\int_0^s K(t)f(t)\,dt,\qquad
G(s)=-\int_0^s K(t)^2\gamma(t)\,dt.
\end{aligned}
\end{equation*}
Put $u=Kv$ and $d\theta/dt=1/K(s)$. The identity $K'=-hK$ eliminates the quadratic velocity term:
\begin{equation*}
\frac{ds}{d\theta}=u,\qquad
\frac{du}{d\theta}=F'(s)u-G'(s).
\end{equation*}
Equivalently, setting $w=u-F(s)$,
\begin{equation}\label{eq:potential}
s_\theta=w+F(s),\qquad w_\theta=-G'(s).
\end{equation}
Here $F(0)=F'(0)=G(0)=G'(0)=0$ and $G''(0)=1$. The coordinate and time changes are real analytic with positive time factor near the equilibrium.

\begin{theorem}\label{thm:analytic}
Let
\begin{equation*}
\xi(s)=s\sqrt{\frac{2G(s)}{s^2}},\qquad \psi=\xi^{-1},
\qquad \sigma(s)=\psi(-\xi(s)),
\end{equation*}
with the square-root branch normalized by $\xi'(0)=1$. The following are equivalent:
\begin{enumerate}[label=(\roman*),itemsep=0.2ex]
\item $(1,0)$ is a center of \eqref{eq:family};
\item $H(\xi):=F(\psi(\xi))$ is even;
\item $F=\Phi\circ G$ for some local real analytic $\Phi$;
\item $F\circ\sigma=F$ and $G\circ\sigma=G$.
\end{enumerate}
The map $\sigma$ is the unique local nontrivial analytic involution fixing $0$, with $\sigma'(0)=-1$, that preserves $G$.
\end{theorem}
The center--evenness equivalence is the classical analytic Li\'enard criterion in the form used by Cherkas \cite{Cherkas1976}. Indeed, the analytic change $X=\xi(s)$ and the positive time factor $d\vartheta/d\theta=\xi'(s)$ transform \eqref{eq:potential} into
\[
X_{\vartheta}=w+H(X),\qquad w_{\vartheta}=-X.
\]
If $H$ is even, this system is reversible under $(X,w)\mapsto(-X,w)$. Conversely the center criterion for this normalized analytic Li\'enard form gives evenness of $H$. Since $G=\xi^2/2$, evenness of $H$ is equivalent to $F=\Phi\circ G$, and to invariance under $\sigma$.

\section{Rational descent and the finite degree bound}\label{sec:rational}
Define the polynomial expressions and analytic auxiliary function
\begin{equation}\label{eq:CDYL}
\begin{aligned}
&C=AB'-BA'-4Bp,\qquad D=csA+Bp,\qquad Y=KB,\\
&L=\frac{C}{AB},\qquad J=\frac{D}{AB}=\frac{cs}{B}+\frac{p}{A}.
\end{aligned}
\end{equation}
A direct derivative calculation gives
\begin{equation}\label{eq:rational-identities}
\frac{Y'}Y=L,\qquad G'=-\frac{4K^2BD}{A},\qquad
F'=Y'+\frac{G'}Y.
\end{equation}
In the nonvanishing layer $C\not\equiv0$, introduce
\begin{equation*}
R=-\frac{4D}{C},\qquad S=\frac{L}{R'}
\quad\text{when $R$ is nonconstant}.
\end{equation*}
The transcendental factor $K$ cancels in these ratios. In particular, $G'=RYY'$ wherever the meromorphic expressions are defined.

\begin{proposition}\label{prop:rational-center}
Outside the layers $C\equiv0$ and $R$ constant, the center criterion of Theorem~\ref{thm:analytic} is equivalent to the existence of a local analytic involution $\sigma$ such that
\begin{equation*}
R\circ\sigma=R,\qquad S\circ\sigma=S,
\qquad\sigma(0)=0,\quad\sigma'(0)=-1.
\end{equation*}
\end{proposition}
\begin{proof}
If the center criterion holds, $Y$ and $Y\circ\sigma$ solve the same analytic initial-value problem $Z'=F'-G'/Z$ with $Z(0)=-1/4$. Thus $Y\circ\sigma=Y$. Differentiating $G\circ\sigma=G$ and $Y\circ\sigma=Y$ shows $R\circ\sigma=R$; the logarithmic derivative gives $L(\sigma)\sigma'=L$, and differentiation of the first invariance identity then gives $S\circ\sigma=S$. Conversely these rational invariances imply $L(\sigma)\sigma'=L$, so $Y\circ\sigma=Y$ by its logarithmic derivative and initial value. Identity \eqref{eq:rational-identities} then gives invariance of $F,G$, proving the center by Theorem~\ref{thm:analytic}.
\end{proof}

The two exceptional layers are essential. If $C\equiv0$, then $Y\equiv-1/4$ and $F=G/Y(0)$. If $R\equiv R_0$, then $R_0\ne0$ and
\[
G=\frac{R_0}{2}\bigl(Y^2-Y(0)^2\bigr),\qquad
F=(1+R_0)\bigl(Y-Y(0)\bigr).
\]
Both layers satisfy the center criterion and belong to condition $\CCfam1$. These cases must be removed before dividing by $C$ or $R'$.

There is also an exact membership test not requiring a complete classification. Write the reduced fractions $R=R_n/R_d$ and $S=S_n/S_d$ and form
\begin{equation*}
\begin{aligned}
J_R(s,z)&=R_n(s)R_d(z)-R_n(z)R_d(s),\\
J_S(s,z)&=S_n(s)S_d(z)-S_n(z)S_d(s).
\end{aligned}
\end{equation*}
Let $Q(s,z)$ be their gcd with every diagonal factor $s-z$ removed; when $S$ is constant take $\gcd(J_R,0)=J_R$. The following theorem gives the resulting finite algebraic test.
\begin{theorem}\label{thm:finite}
A nonexceptional parameter vector gives a center if and only if $Q$ is nonconstant and its lowest-degree nonzero homogeneous part $Q_m$ satisfies
\begin{equation}\label{eq:Qcriterion}
Q_m(1,-1)=0.
\end{equation}
\end{theorem}
\begin{proof}
At a center the graph $z=\sigma(s)$ is a non-diagonal analytic branch of both equality curves, with tangent $z=-s$. For the converse, choose a real rational generator $W$ of $\mathbb R(R,S)$ (by L\"uroth over $\mathbb R$), and replace $W$ by $1/W$ if it has a pole at zero. The local normal-form theorem for an analytic map of multiplicity $e$ provides an analytic coordinate $\nu$ with $\nu(0)=0$, $\nu'(0)\ne0$, in which $W-W(0)=a\nu^e$, $a\ne0$. Locally, the common equality-curve germ of $R$ and $S$ is that of $W(s)=W(z)$: the field-generation identities express $W$ rationally in $R,S$, while the reverse inclusions are automatic on the locus where the rational functions are defined. Its distinct analytic branches are therefore
\[
 z=\nu^{-1}(\zeta\nu(s)),\qquad \zeta^e=1,
\]
whose tangent lines are $z=\zeta s$. The diagonal $\zeta=1$ is removed in $Q$. Consequently $Q_m(1,-1)=0$ is equivalent to the occurrence of the factor $z+s$ in the lowest-degree tangent cone, and hence to $e$ even. The generator and its local normalization can be chosen real and the branch $z=\nu^{-1}(-\nu(s))$ is a real analytic involution with derivative $-1$; it preserves $W$, and thus $R,S$. Reducing denominators only removes base factors of their rational presentations, not this analytic branch. Proposition~\ref{prop:rational-center} gives the center.
\end{proof}

\begin{theorem}\label{thm:degree}
For any nonexceptional center, the subfield
\begin{equation*}
\calK=\CC(R,S)\subset\CC(s)
\end{equation*}
is $\CC(W)$ for a rational function $W$, and
\begin{equation*}
d=[\CC(s):\CC(W)]=\deg W,\qquad 2\le d\le5.
\end{equation*}
If $d=2$ the induced deck involution is M\"obius, whereas for $d\ge3$ one has $\CC(W)=\CC(R)$.
\end{theorem}
\begin{proof}
L\"uroth's theorem gives a rational generator $W$. The local negative-slope involution fixes $W$, so its local degree at zero is even, in particular $d\ge2$. Writing $A_6=-a_5$ and $B_2=c$, their highest coefficients satisfy $[s^7]C=(-2c+6c-4c)a_5=0$ and $[s^7]D=(-c+c)a_5=0$. Also $C(0)=D(0)=0$, so $s$ divides both. Thus $\deg R\le5$ after reduction. Since $R\in\CC(W)$, $\deg R=d\deg(r)$ for $R=r\circ W$, proving the upper bound. For $d\ge3$, the inequality $d\deg(r)\le5$ forces $\deg(r)=1$, so $R$ itself is a generator. In degree two, the extension is Galois; its nontrivial automorphism belongs to $\operatorname{PGL}_2(\CC)$ and has the normalized form $\sigma(s)=-s/(1+bs)$.
\end{proof}

A further local ramification observation sharpens the argument: for a nonexceptional center the local degree of the common generator at zero is exactly two. Writing $m=\ord_0 C\ge1$ and using $\ord_0D=1$, one finds that for $m>1$ the orders of $R$ and $S$ are $1-m$ and $2m$, while for $m=1$, if $j=\ord_0(R-R(0))$, then $\ord_0S=2-j$. The common even local degree divides two in either case. This prevents a hidden local-degree-four branch in later ramification accounting.

\paragraph*{Cancellation and reduced rational invariants}
The descent requires cancellation of the integrating factor in
\eqref{eq:rational-identities}, reduction of the rational functions to
coprime numerator--denominator form, and clearing denominators before
forming the bivariate greatest common divisor in the finite membership
test. In particular, the non-diagonal branch with tangent $-1$ must
survive the reduced gcd. The degree estimate also requires the
highest-degree cancellation in $C,D$ and the common-factor cancellation
before applying $\deg R\le5$.

\section{Residues, fibers and ramification}\label{sec:residue}
Set
\begin{equation}\label{eq:T}
T=(1-R)L
=\frac{3B'+2g}{B}-\frac{A'}A,
\qquad\omega=T(s)\,ds.
\end{equation}
For $d\ge3$ the equality $\CC(R,S)=\CC(R)$ implies $T/R'=(1-R)S\in\CC(R)$. Consequently there exists a rational $\eta(w)$ satisfying
\begin{equation*}
\omega=\eta(R(s))\,dR(s)=R^*(\eta(w)\,dw).
\end{equation*}
This simple differential is the link between field descent and finite combinatorics. In the nondegenerate real stratum $cg\ne0$, $\Delta=g^2+c>0$, the two roots of $B$ are distinct. Put
\[
\beta_\pm=\frac{g\mp\sqrt\Delta}{2c},\qquad
\tau=\frac{2g}{\sqrt\Delta},
\]
with the signed root convention fixed before any residue comparison. At a root $\alpha$ of $A$ disjoint from $B$, of multiplicity $m_\alpha$, the residue is $-m_\alpha$. At $\beta_\pm$, if $m_\pm=\ord_{\beta_\pm}A$, one obtains
\begin{equation}\label{eq:residues}
\Res_{\beta_\pm}\omega=3\mp\tau-m_\pm.
\end{equation}
An effective zero multiplicity $6-\deg A$ at infinity must enter the global accounting whenever $\deg A<6$. In particular, the total ordinary negative mass (including the effective infinity multiplicity and excluding multiplicities absorbed at roots of $B$) is
\begin{equation}\label{eq:mass}
K_{-}=6-m_+-m_-,\qquad N_{-}\le K_{-}\le6.
\end{equation}

\begin{lemma}\label{lem:res-pull}
Let $R$ have local degree $e_p$ at $p$ and write $w_0=R(p)$. Then
\begin{equation*}
\Res_p(R^*\eta\,dw)=e_p\Res_{w_0}(\eta\,dw).
\end{equation*}
Moreover $\ord_p(R^*\omega_0)=e_p(\ord_{w_0}\omega_0+1)-1$ for a meromorphic differential $\omega_0$ on the target.
\end{lemma}
\begin{proof}
Use local coordinates in which $w-w_0=t^{e_p}$ up to a nonvanishing unit. The logarithmic term $\rho\,dw/(w-w_0)$ pulls back to $e_p\rho\,dt/t$; the general order formula follows by pulling back $w^m\,dw$.
\end{proof}

The fiber identity and Riemann--Hurwitz formula are respectively
\begin{equation*}
\sum_{R(p)=w_0}e_p=d,\qquad
\sum_{p\in\PP^1}(e_p-1)=2d-2.
\end{equation*}
Since the target differential has only simple poles in the relevant nonexceptional classification, every point in the preimage of a nonzero-residue target pole is itself a nonzero-residue source pole. Thus source poles in the same fiber have residue ratios equal to the positive integers $e_p$. For example, in degree four the only fiber partitions are
\[
[1,1,1,1],\quad [2,1,1],\quad [2,2],\quad[3,1],\quad[4],
\]
which imply the respective nonzero residue ratios $1:1:1:1$, $2:1:1$, $1:1$, $3:1$, and one fully ramified pole.

\begin{proposition}\label{prop:mass}
For nonexceptional degrees $3\le d\le5$, the positive-residue source roots are contained in one target fiber, either together or with only one surviving positive pole. The number $r_-$ of negative-residue target fibers is one or two. If the two roots of $B$ lie in the positive fiber, then
\begin{equation}\label{eq:two-positive}
N_-+d(1-r_-)\ge1.
\end{equation}
If just one root contributes a positive pole, set $\varepsilon=1$ when the other root is a negative pole, and $\varepsilon=0$ when its residue cancels. Then
\begin{equation}\label{eq:one-positive}
N_-+\varepsilon-1+d(1-r_-)\ge1.
\end{equation}
\end{proposition}
\begin{proof}
With two positive roots in one fiber, their fiber uses $d-2$ units of total ramification. The $r_-$ negative fibers, supported at $N_-$ distinct source points, use $r_-d-N_-$. The obligatory quadratic local ramification at zero uses one more unit. The total is bounded by $2d-2$, giving \eqref{eq:two-positive}. The singleton-positive case contributes $d-1$ units and yields \eqref{eq:one-positive}. Combining these inequalities with \eqref{eq:mass} gives $r_-\in\{1,2\}$.
\end{proof}

The polynomiality reconstruction must accompany residue matching: if a rational configuration is formally admissible, it need not yet correspond to a quintic vector field. We use the exact recovery identity
\begin{equation*}
p=-\frac{A}{4}\left(\frac{R}{1-R}T+\frac{4cs}{B}\right),
\qquad q=A+sp,
\end{equation*}
and impose $\deg p,\deg q\le5$ with the prescribed coefficients at $s=0$. This is the essential bridge between abstract fiber configurations and ten-parameter polynomial systems.

\paragraph*{Compatibility with polynomial reconstruction}
For any proposed rational generator, the poles of $\omega$ must be
accounted for at every finite root of $A$ and $B$ and at infinity,
including cancellations and multiple roots. The possible fiber
partitions must satisfy both the pullback-residue multiplicities and
the Riemann--Hurwitz budget $\sum_p(e_p-1)=2d-2$.
These conditions make the list of integer ramification patterns finite,
while polynomial reconstruction imposes an additional realizability condition, treated in Section~\ref{sec:classification}.

\section{Exhaustive parameter classification and completeness}\label{sec:classification}
We now apply the criterion of Sections~\ref{sec:lienard}--\ref{sec:rational} and the residue restrictions of Section~\ref{sec:residue} to the whole parameter space. Here the \emph{degree} of a nonexceptional family means
\[
 d=[\CC(s):\CC(W)]=\deg W,
\]
as in Theorem~\ref{thm:degree}; it refers to the common rational generator, not to the degree of the original vector field, which is five throughout. This integer $d$ is distinct from the free shape parameter denoted by $d$ in conditions $\CCfam{23}$, $\CCfam{24}$ and $\CCfam{36}$. The four elementary conditions $\CCfam1$--$\CCfam4$ are treated first. The remaining conditions are organized into nine degree-two families ($\CCfam5$--$\CCfam{13}$), sixteen degree-three families ($\CCfam{14}$--$\CCfam{29}$), nine degree-four families ($\CCfam{30}$--$\CCfam{38}$), and four degree-five families ($\CCfam{39}$--$\CCfam{42}$). These labels describe the reconstruction branches; overlaps and lower-degree specializations remain included. We first specify the order of the branch analysis and then give its coefficient and residue calculations.

\paragraph*{The calculation order and the destination of every structural branch}\label{sec:necessity-order}
Start from an arbitrary real center, put $M=4a_1+2a_7$, and compute $c=a_1+1$, $g=a_6$, $p,A,B,C,D$ from Sections~\ref{sec:lienard}--\ref{sec:rational}. We use the following ordered decision procedure. It is not assumed at the start that $c,g,M+2,M+6$ are nonzero.
\begin{enumerate}
\item Test $C=(M+2)D$ by its five polynomial coefficients. If it holds, the parameter vector belongs to $\CCfam{1}$; this includes $C=0$ and both exceptional exponent layers.
\item On the complement, first separate $g=0$. If $c\ne0$, the logarithmic-invariant argument below leaves only A, hence $\CCfam{2}$. If $c=0$, the polynomial-composition argument additionally leaves G, hence $\CCfam{3}$.
\item If $c=0,g\ne0$, divide only after proving $B\mid A$ and $B\mid p$. The reduced constant-$B$ problem gives the reversible lift H, hence $\CCfam{4}$, or returns to F.
\item Now $cg\ne0$. Split $\Delta=g^2+c$ into its three signs. The nonpositive cases return to F or contradict the pole-order and ramification bounds. Only then normalize the two distinct real roots and enumerate the integer residues.
\item For a quadratic full invariant field, solve all fixed rational covariance systems; the nine surviving source families give $\CCfam{5}$--$\CCfam{13}$. For degrees three, four and five, use the finite mass partitions below and reconstruct all original coefficients. Their retained families give $\CCfam{14}$--$\CCfam{29}$, $\CCfam{30}$--$\CCfam{38}$ and $\CCfam{39}$--$\CCfam{42}$, respectively, after the stated inclusions and boundary mergers.
\end{enumerate}
The arrows describe the first applicable stratum. Some boundary points also lie in earlier families. A source label is not automatically a new center condition: \ref{app:records} gives every one of the 69 source labels and its final condition number.

\paragraph*{Coefficient comparison in original parameters, before any division}\label{sec:F-recursion}
To recover condition $\CCfam{1}$ without dividing by a parameter, set $M=4a_1+2a_7$ and write
\begin{center}
\begin{adjustbox}{max width=\linewidth}
$\displaystyle A=\sum_{j=0}^6A_js^j,\quad p=\sum_{j=1}^5p_js^j,\quad
A_0=1/4,\ A_1=g,\ A_2=a_7-a_1,\ p_1=c-1,\ p_5=-A_6.$
\end{adjustbox}
\end{center}
An index outside these ranges denotes zero. Let $E_j$ be the coefficient of $s^j$ in
\[
\mathcal E=B[A'+(M+6)p]+A(g+Mcs).
\]
The exact convolution rule is
\begin{equation*}
\begin{aligned}
E_j={}&-\frac{j+1}{4}A_{j+1}-g(j-1)A_j+c(M+j-1)A_{j-1}\\
&+(M+6)\left(-\frac{p_j}{4}-gp_{j-1}+cp_{j-2}\right).
\end{aligned}
\end{equation*}
The $j=0,1$ coefficients vanish by normalization, and the top coefficient vanishes by $p_5=-A_6$. For $j=2,\ldots,6$ the convolution gives $E_j=R_j$, with $R_j$ defined in Theorem~\ref{thm:main}. Thus
\[
B[A'+(M+6)p]+A(g+Mcs)=\sum_{k=2}^6R_ks^k,
\]
and the five remaining coefficient equations are exactly $R_2=\cdots=R_6=0$. In the expanded notation used below, $f_1=-4R_2$, $f_j=4R_{j+1}$ for $2\le j\le4$, and $f_5=R_6$. This identifies the compact condition in the theorem with the full coefficient calculation. The same identity supplies the integrating factor and first integrals in the sufficiency proof.

For example, $E_2=0$ has the nonzero constant coefficient $-3/4$ on $A_3$, so solve for $A_3$. Next, $E_3,E_4$ solve for $A_4,A_5$ with constant pivots. The last two equations must be retained: the coefficient of $A_6$ in $E_5$ is $M/4$, which may vanish. Dividing by $M$ would therefore omit a parameter stratum. The equations as displayed also retain $M=-2,-6$ without exceptional denominator clearing.

\paragraph*{The constant and vanishing layers: condition F}
Keep $C,D,L,R$ from \eqref{eq:CDYL} and write $M=4a_1+2a_7$. Direct expansion gives $D'(0)=1/4$ and $C'(0)=(M+2)/4$. Thus a constant quotient $R=-4D/C$, or $C\equiv0$, is equivalent to
\begin{equation*}
C=(M+2)D.
\end{equation*}
The coefficient calculation is explicit:
\[
C-(M+2)D=\frac14f_1s^2-\frac14f_2s^3-\frac14f_3s^4-\frac14f_4s^5-f_5s^6.
\]
Thus its five coefficient equations are precisely $\CCfam{1}$. Equivalently,
\begin{equation}\label{eq:Fidentity}
B\bigl((M+6)p+A'\bigr)+A(Mcs+g)=0.
\end{equation}
At $M=-2$ the first identity reads $C=0$ and is retained. At $M=-6$ the second identity loses its $p$ term and is retained as well. Neither layer is divided out. All remaining cases below assume a nonconstant $R$.

\paragraph*{Why only degrees two through five can occur}
The involution preserves $R$ and $S=L/R'$. By Lüroth's theorem their common rational function field is $\C(R,S)=\C(W)$; write $d=[\C(s):\C(W)]$. Choose $W$ regular at zero. Its local degree there is even, since $W\circ\sigma=W$ and $\sigma'(0)=-1$. The potential seventh-degree terms in $C,D$ cancel, and $s$ divides both. Consequently $\deg R\le5$, and
\begin{equation*}
2\le d\le5.
\end{equation*}
For $d\ge3$, the equality $\deg R=d\deg(R\text{ as a function of }W)$ implies that $R$ itself generates the common field. For $d=2$ there is a real Möbius deck involution $\sigma(s)=-s/(1+bs)$.

There is no unexamined local-degree-four stratum. To see this directly, let $m=\operatorname{ord}_0 C\ge1$ and use $\operatorname{ord}_0D=1$. If $m>1$, then $\operatorname{ord}_0R=1-m$ and $\operatorname{ord}_0S=2m$; the even local degree $e$ of $W$ divides both, hence divides two. If $m=1$, write $j=\operatorname{ord}_0(R-R(0))>0$. Then $\operatorname{ord}_0S=2-j$, so $e$ divides $j$ and $2-j$, and again $e=2$. Every non-F center therefore has exactly quadratic local ramification at zero.

The rational differential organizing the remaining classification is $T(s)\,ds=(1-R)L\,ds$, as defined in \eqref{eq:T}.
For $d\ge3$ this descends to $\eta(R)\,dR$. Keep the fixed roots $\beta_\pm=(g\mp\sqrt\Delta)/(2c)$ and the \emph{signed} quantity $\tau_0=2g/\sqrt\Delta$ from Section~\ref{sec:residue}, with $\sqrt\Delta>0$. Their residues are exactly $3\mp\tau_0-m_\pm$ by \eqref{eq:residues}. For the combinatorial enumeration alone, relabel the two roots so that $\widehat\beta_+$ is the root whose coefficient is $3+|\tau_0|$, and $\widehat\beta_-$ is the other; put $\widehat m_\pm=\ord_{\widehat\beta_\pm}A$ and $\widehat\tau=|\tau_0|>0$. Thus the two residues in this \emph{relabelled} order are $3+\widehat\tau-\widehat m_+$ and $3-\widehat\tau-\widehat m_-$. In the integer records below we abbreviate the hatted indexing variables as $\tau,m_+,m_-$; this is only a choice of labels, not a change to \eqref{eq:residues}. All other nonzero residues are negative integers. The total ordinary negative mass is $K=6-m_+-m_-\le6$, and its support size is $N\le K$. The effective multiplicity at infinity is $6-\deg A$ and must be included.

\paragraph*{The discriminant and zero-parameter layers}
\paragraph*{$\Delta=0$, $g\ne0$}
The differential \eqref{eq:T} has one double pole and all other poles are at most simple. A target pole of order $j\ge2$ pulls back at local degree $e$ to order $e(j-1)+1$. A unique source pole would have to be a singleton fiber with $e=d\ge2$, giving order at least three. This contradicts the double pole. Thus only F remains.

\paragraph*{$\Delta<0$, $g\ne0$}
The two conjugate roots of $B$ have residues with nonzero imaginary parts. They cannot share a fiber with ordinary negative integer residues, and their residues cannot be positive rational multiples of each other. Each must therefore be a singleton fully ramified fiber. These consume the entire Riemann--Hurwitz budget $2d-2$, leaving no ramification at zero. The degree-two involution gives the same obstruction. Thus the non-F real center locus with $g\ne0$ has $\Delta>0$. The case $g=0$, including its negative-discriminant part, is treated separately next; reversible family A is not excluded by this argument.

\paragraph*{$g=0$, $c\ne0$: conditions F and A}
Now $T\,ds=d\log\Psi$, with $\Psi=B^3/A$. In degree two, the involution $-s$ gives $A$ even and $p$ odd, exactly A. A nonlinear involution can only fix one root $\beta$ of $B$ as a zero of $\Psi$. With $z=s/\beta$, $r=\beta^2$, the exceptional shape is
\begin{center}
\begin{adjustbox}{max width=\linewidth}
$\displaystyle A=\tfrac14(1-z)(1+z)^3(1-2z+dz^2),\qquad
\beta p(\beta z)=(1/4-r)z+j_2z^2+j_3z^3+j_4z^4+dz^5/4.$
\end{adjustbox}
\end{center}
The three vanishing pole coefficients at $z=-1$ form a linear system of determinant $96$, giving
\[
j_2=-3r,\qquad j_3=d/4-3r-3/4,\qquad j_4=d/2-r-1/2.
\]
Substitution makes $R=8r/(8r+1-d)$ constant, or $C=0$ when the denominator vanishes, so this adds only F.

For degrees three to five, logarithmic descent gives $\Psi^n=H(R)$ for some positive integer $n$. The associated cyclic rational extension has degree $k$, with $d=kf$. Since $\Psi(0)=-1/16$ is neither zero nor infinity, a nontrivial power map is unramified at the source value over zero; the required ramification gives $f\ge2$. Thus degrees three and five have $k=1$, while the only additional degree-four possibility is $k=f=2$, already in degree two. The bound $\deg\Psi\le6$ then makes $\Psi$ Möbius-equivalent to $R$; the possible degree-six quadratic outer map in degree three would use two triple fibers and leave no ramification at zero.

Let $a$ be the value of $\Psi$ corresponding to $R=1$. If $a=0$ or $\infty$, polynomiality forces equal multiplicities at the two roots of $B$, hence $A$ even and $p$ odd. For a finite nonzero $a$, the only cancellation needing attention has $A=(s-\beta)^3E$, $\deg E\le3$, with $E'(0)/E(0)=3/\beta$. A corresponding finite-fiber polynomial $Q=(s+\beta)^3-aE$ has $Q'(0)/Q(0)=3/\beta$; if supported only at $\beta$ it would have the incompatible ratio $-j/\beta$. The fiber is therefore at infinity. Writing $J=\gcd(A,B^3)$ gives $B^3-aA=\lambda J$, $\deg J=6-d$, $J'(0)=0$. Equal root multiplicities exclude odd $6-d$ and give A in degree four. Hence this entire zero-$g$ layer contains only F and A.

\paragraph*{$c=g=0$: the additional family G}
Here $B=-1/4$ and $T=-A'/A$. A nonlinear degree-two involution would create a finite pole in the polynomial identity $A\circ\sigma=A$, so $A$ must be constant, giving F. The involution $-s$ gives A. For $d\ge3$, write the reduced quotient $R=u/v$; then $A'=(u-v)H$, $p=-uH/4$, with $\deg H\le5-d$. The conditions $p(0)=0$, $p'(0)=-1$ and even ramification force $H$ to have a simple zero, so $d\le4$. The positive pole at infinity makes the generator polynomial. Degree three yields
\[
W=s^2+ts^3,\qquad A=1/4+\rho W-3\kappa W^2,\qquad
p=(-1/2+\kappa W)W',
\]
which is G after coefficient comparison. The nontrivial degree-four cyclic case gives an even quadratic inner generator and hence A; the remaining degree-four case would require $\deg p=7$. Thus $g=0$ gives exactly F, A and G.

\paragraph*{$c=0$, $g\ne0$: the additional family H}
The differential is $d\log(B/A)$. If $B\nmid A$, the simple finite pole of $A/B$ and its infinity fiber force $A/B$ to generate the common polynomial field. Polynomiality of $p$ would make $p$ bounded at infinity, contradicting $p(0)=0$, $p'(0)=-1$. Thus $B\mid A$. A similar singleton-fiber argument gives $B\mid p$. The original vector field has the factor $x+4gy$, an analytic unit equal to one at $(1,0)$. Dividing by it gives
\[
x'=-x^3y+b_2x^2y^2+b_3xy^3+b_4y^4,\qquad
y'=-1/4+x^4/4+d_2x^2y^2+d_3xy^3+d_4y^4.
\]
Multiplication by the positive unit $x$ returns to the already classified $c=g=0$ quintic layer, now with $a_5=a_{10}=0$. Its center conditions are
\begin{center}
\begin{adjustbox}{max width=\linewidth}
$\displaystyle b_2=b_4=d_3=0,
\quad\text{or}\quad
3d_3=(1-2d_2)b_2,\quad2d_4=(1-d_2)b_3,\quad(2d_2-3)b_4=0.$
\end{adjustbox}
\end{center}
The first alternative restores to H. The second has integrating factor $x^{2d_2-3}/(x+4gy)$ in the original system and belongs to F. Thus this layer gives exactly F and H.

\paragraph*{The degree-two common field: all nine rational families}
Assume $cg\ne0$ and a degree-two common field. Its deck involution is $\sigma=-s/(1+bs)$. If it swaps the roots of $B$, equal residues give $n=m_+-m_-=2\tau\in\{1,\ldots,6\}$, $b=4g$ and $c/g^2=16/n^2-1$. Put $z=gs$, $X=1+2z$, $Q=1+4z$, $l_+=1-2(4-n)z/n$, $l_-=1+2(4+n)z/n$. Every admissible $A$ has the form
\[
A=l_+^n\left\{\tfrac14X^{6-n}+\sum_{j=1}^{\lfloor(6-n)/2\rfloor}u_jX^{6-n-2j}z^{2j}\right\}.
\]
If the roots are not swapped, their residues force an integer $\tau=j$. The full candidate set is $j=-1,1,3,4,5,6,7,8,9$; the zero-$c$ value $j=2$ has already been treated. The upper bound follows from the effective multiplicity constraint $m_++2m_-+j-3\le6$, also when the image point is at infinity.

Set $\rho=1/g^2$, $\varkappa=c/g^2$, $v(z)=p(z/g)/g$, and
\[
A=1/4+z+\sum_{j=2}^6A_jz^j,\qquad
v=(\varkappa-\rho)z+\sum_{j=2}^5v_jz^j.
\]
The exact identities are $A(\sigma)/A=F_B(\sigma)/F_B$ and $Z(\sigma)\sigma'=Z$, with $Z=\varkappa z/B+v/A$ and $A_6+v_5=0$. Here $F_B=l_+^{3+j}l_-^{3-j}$ and $l_+=1-2(2-j)z/j$, $l_-=1+2(2+j)z/j$ in the fixed-root case. After clearing only powers of the position denominator, these are linear equations with constant rational matrices for $A_j,v_j,\rho$. Their exact row reduction is exhaustive; it has no parameter-dependent rank assumption.

For swapped roots, $n=1,2,3$ give $\CCfam{5},\CCfam{6},\CCfam{7}$; $n=4$ forces $\rho=0$, $n=5$ forces $\rho=-4/25$, and $n=6$ is inconsistent. The latter alternatives give no additional finite real $g$. For fixed roots, $j=-1$ gives F, $j=1,3,4,5,6,7$ give $\CCfam{8}$ through $\CCfam{13}$, $j=8$ forces $\rho=-3/4$, and $j=9$ is inconsistent. Exact left-row combinations certify the inconsistent systems; substitution checks every retained solution and $A_6+v_5=0$. The result is exactly nine non-F real degree-two families.

\paragraph*{The shared mass and ramification bounds for degrees three to five}
The target differential $\eta(R)dR$ has only simple poles; a higher-order target pole would pull back to a forbidden higher-order source pole. Positive fibers contain only roots of $B$. Two different singleton positive fibers would use all $2d-2$ ramification, so there is one positive fiber: either both roots lie in it, or one root is fully ramified. Let $r_-$ be the number of negative target fibers.

When both roots lie in the positive fiber with local degrees $j,d-j$, its ramification is $d-2$. The negative fibers contribute $r_-d-N$; keeping ramification at zero gives
\[
N+d(1-r_-)\ge1.
\]
For a singleton positive fiber, let $\epsilon=1$ if the other root is a negative pole and $\epsilon=0$ if its residue cancels. The corresponding inequality is
\[
N+\epsilon-1+d(1-r_-)\ge1.
\]
Together with $N\le K\le6$ these leave only $r_-=1,2$. The target residue times each local degree must also be an integer multiplicity in $A$; support counts alone are not sufficient.

Polynomiality is imposed using the exact identity
\begin{equation*}
p=-\frac A4\left(\frac{R}{1-R}T+\frac{4cs}{B}\right),\qquad q=A+sp.
\end{equation*}
If $R=1$ is a target pole, its source local degree $e$ and $A$ multiplicity $a$ require $a\ge e+1$, including the effective infinity multiplicity $a=6-\deg A$. If $R=1$ were a regular nonzero value of $\eta$, an infinity point would force a nonzero $s^6$ coefficient of $q$, while a finite point could only be a canceled root of $B$ and would use a second fully ramified fiber. Both are impossible. Thus $R=1$ lies at a zero or pole of $\eta$. This treats canceled-residue cases, including $\tau=\pm3$; it does not assume that $L$ itself has only simple poles.

\paragraph*{Reconstruction of a source record: the exact equations that are solved}\label{sec:reconstruction-steps}
A partition determines a rational generator $W=N/V$ and the divisor of $A$. It does not yet determine a center condition in the ten original coefficients. The reconstruction has the following steps.
\begin{enumerate}
\item Choose the root coordinate $z=s/\beta$, $\beta\ne0$, and retain $r=\beta^2>0$. Normalize $N(0)=V(0)=1$ and impose $N'(0)=V'(0)$, which is $W'(0)=0$. For a linear ordinary factor this fixes its linear coefficient; for a quadratic factor the remaining coefficient stays free. Test $W_{zz}(0)\ne0$ separately.
\item Construct $A(z)$ from the assigned root multiplicities and the remaining target fiber. For a finite type-I record this reads
\[
A(z)=\frac{(1-z)^m(1-z/b)^nP(z)^k[N(z)-hV(z)]}{4(1-h)}.
\]
The power of each factor is its exact source multiplicity. This makes the mass equation and $A(0)=1/4$ explicit. The repeated-fiber cases use the corresponding power or quadratic outer factor; no simple-root assumption is made there.
\item \emph{Distinguish the two descended differentials.} Write
\[
\begin{aligned}
T(s)&=\vartheta(W(s))W'(s),&
J(s)&=\jmath(W(s))W'(s),\\
R(s)&=\mathcal R(W(s)).
\end{aligned}
\]
where $J$ is the coefficient-reconstruction differential defined in \eqref{eq:CDYL}. The identities $J=-RL/4$ and $T=(1-R)L$ yield the exact rational relation
\begin{equation}\label{eq:two-descents}
\jmath(w)=-\frac{\mathcal R(w)}{4(1-\mathcal R(w))}\,\vartheta(w).
\end{equation}
In particular, the rational function $\vartheta$ used to match residues of $T\,ds$ is \emph{not} the function $\jmath$ used to reconstruct the coefficients. Write the former in its allowed partial-fraction form
\[
\vartheta(w)=\sum_{\ell}\sum_{j=1}^{q_\ell}
\frac{b_{\ell j}}{(w-w_\ell)^j}+
\sum_{j=0}^{q_\infty}d_jw^j.
\]
The pole-order constraints determine all upper limits, and an absent polynomial part is interpreted as zero. Match the residues of $T\,ds$ at all relevant source poles. The outer rational map $\mathcal R$ and the relation \eqref{eq:two-descents} then determine $\jmath$. Since $J'(0)=-4$, with $z=s/\beta$ and $r=\beta^2$ one must impose
\[
\jmath(1)W_{zz}(0)=-4r.
\]
After the outer map and pole locations are fixed, the partial-fraction residue matching is linear in its coefficients; any remaining parameter and denominator branches must still be treated explicitly.
\item In the source coordinate put $\widehat B(z)=B(\beta z)$, $P(z)=\beta p(\beta z)$, and $Q(z)=(A+sp)(\beta z)$. The coefficient-reconstruction identity (with $A=A(\beta z)$) is
\begin{equation}\label{eq:reconstruction-J}
P(z)=A(\beta z)\left[\jmath(W(z))W_z(z)
          -\frac{crz}{\widehat B(z)}\right],
\qquad Q(z)=A(\beta z)+zP(z).
\end{equation}
Verify that all finite poles cancel, $\deg P,\deg Q\le5$, $P_1=(c-1)r$, and $A_6+P_5=0$. Infinity must be included in this check, not handled only by a finite-root test.
\item Finally restore every original coefficient:
\[
a_j=\frac{[z^j]P}{\beta^{j+1}},\qquad a_{j+5}=\frac{[z^j]Q}{\beta^j},\quad1\le j\le5.
\]
Equality of these ten rational functions is the test for a root-label duplicate, a family containment or a boundary merger. A common discriminant value alone is not such a test.
\end{enumerate}
This procedure also explains completeness of the coefficient calculation. Every center supplies $W,A,\mathcal R,\vartheta,\jmath$ satisfying the listed mass, pole and normalization equations. Conversely the elimination solves the whole finite linear system with its free variables, not selected numerical points. Any extra solution would therefore require an omitted partition or a missing parameter-denominator stratum; both are treated in the branch analysis below.

\paragraph*{A worked degree-three record and its final number}
For D3-I-02 the indices are $(k,e,m,n)=(1,1,1,1)$, $b=1/3$, $\tau=4$. The actual polynomials are
\[
N=(1-3z)^2(1+3z),\quad V=(1-z)^3,\quad W=N/V,\quad W_{zz}(0)=-24,
\]
\[
A=\frac{(1-z)(1-3z)(1+3z)(N-hV)}{4(1-h)}.
\]
Here $A_z(0)=-1$, so $g\beta=-1$. Solving the allowed target-differential equations leaves a free coefficient $t$ and gives, among the restored source coefficients,
\begin{center}
\begin{adjustbox}{max width=\linewidth}
$\displaystyle\begin{aligned}
P_1&=-r-3/4,\quad P_2=3r,\quad Q_1=-1,\\
P_3&=-\frac{3(2hr-3h-2r+24t+9)}{2(h-1)},\qquad
Q_2=-\frac{2hr+3h-2r-9}{2(h-1)}.
\end{aligned}$
\end{adjustbox}
\end{center}
Set $k_0=1/(h-1)$, $u_0=t/(h-1)$ and $\beta=-1/g$. Coefficient restoration gives
\begin{center}
\begin{adjustbox}{max width=\linewidth}
$\displaystyle\begin{aligned}
a_1&=-1-3g^2/4,\quad a_2=-3g,\quad a_6=g,\\
a_3&=-3g^2+(9/2-9k_0-36u_0)g^4,\qquad
a_7=-1+(3k_0-3/2)g^2.
\end{aligned}$
\end{adjustbox}
\end{center}
Here $k_0,u_0$ are the parameters denoted $k,u$ in condition~14. The remaining coefficients follow from the same $P,Q$ coefficient equations and are precisely the displayed map of $\CCfam{14}$. Its $k_0=0$ boundary includes D3-B-01. Thus this example exhibits the chain from a residue partition, through a fixed linear system and polynomiality, to an actual numbered original-parameter condition.

\paragraph*{Degree three: all integer possibilities and resulting conditions}
If both roots lie in the positive fiber and $r_-=2$, then $N\ge4$ and $K\le5$. The ordinary cubic partitions are $[3]$, $[2,1]$ and $[1,1,1]$. The only admissible pair of negative fibers is one singleton $[3]$ of multiplicity $k=1,2$ and one three-point unramified fiber of multiplicity one at each point. With $j=1,2$ and $m=0,\ldots,3-k$, the residue equation is
\[
\tau=j(1+k/3)-3+m.
\]
This gives ten type-III source records. With $r_-=1$, negative $R=1$ gives their $h=0$ limits; positive $R=1$ leaves two infinity limits from $k=1$, $m=n=1$.

For a singleton positive root and $r_-=2$, the possible ordinary support sizes of the two negative fibers are $(3,1)$, $(3,2)$ and a collision $(2,2)$. The first gives type I with
\[
k=1,2,3,\qquad e=1,2,\qquad m+n=3-k,\qquad \tau=3k/e+m,
\]
namely twelve source records. The other support sizes have $K=5$, $m+n=1$ and give two type-II records with $\tau=3,4$ and a free quadratic factor. All other partitions violate either the support bound or integer mass six. With $r_-=1$, the negative fiber gives $h=0$ limits; positive $R=1$ gives six type-I infinity boundaries from $(k,m,n)=(1,1,1),(1,2,0),(2,1,0)$, $e=1,2$, and one type-II boundary.

Thus degree three has 24 finite source records and nine boundary records. Exact coefficient identifications, inclusions and boundary mergers give conditions $\CCfam{14}$--$\CCfam{29}$. The record-by-record assignment is displayed in \ref{app:records}.

\paragraph*{Degree four: all integer possibilities and resulting conditions}
For two roots in the positive fiber, $r_-=2$ requires $N\ge5$, $K\le5$. A two-point partition $[3,1]$ has minimum mass four, $[2,2]$ mass two, and a three-point partition $[2,1,1]$ mass four. The only permissible pair of negative fibers is therefore a singleton $[4]$ of multiplicity one and a four-point unramified fiber. With $j=1,2,3$, $m=0,1$ and $\tau=5j/4-3+m$, there are six type-III records. The $r_-=1$ negative limit gives their $h=0$ boundary, and the positive limit would require mass at least six at the two roots, leaving no ordinary negative mass.

For one positive root, $r_-=2$ gives $\epsilon=1$, $N\ge5$. Support sizes $(4,1)$ give type I with $e=1,2,3$, $k=1,2$, $m+n=2-k$, $\tau=1+4k/e+m$, or nine records. Support sizes $(4,2)$ or a collision $(3,2)$ give one type-II record with mixed partition $[B:2,1,1]$ and mass six. The remaining $(2,3)$ partition factors through a quadratic inner generator and is already in degree two. Its possible local-degree degeneracy cannot hide a center: for the type-II shape at $d=11/3$, $W'''(0)=-512/9\ne0$, so the local degree would be odd.

With $r_-=1$, a negative limit requires $(m_+,m_-,k)=(0,5-e,e+1)$ and gives finite type-I $h=0$ limits; the positive limit gives $(5,0,1)$, $e=1,2,3$, or three infinity boundaries. Thus degree four has sixteen finite records and three boundary records. Removing whole-family degree-two inclusions and duplicate root labels leaves conditions $\CCfam{30}$--$\CCfam{38}$.

\paragraph*{Degree five: all four remaining conditions}
Two roots in the positive fiber would require either $N\ge6$ with zero mass at both roots, contradicting polynomiality, or a negative singleton limit with $K\ge5+N$, again incompatible. A single positive root with $r_-=1$ also fails the mass bound. For $r_-=2$ the only possibility is
\[
\epsilon=1,\qquad N=K=6,\qquad m_+=m_-=0.
\]
The ordinary negative residues are $-1$. A quintic ordinary fiber has equal local degrees dividing five; its support size is one or five. The total support six forces five unramified points and a mixed fiber $[B:5-e,\alpha:e]$, $e=1,2,3,4$. These are precisely $\CCfam{39}$--$\CCfam{42}$.

Their direct reconstruction, before returning to $g$, uses $z=s/\beta$:
\begin{center}
\begin{adjustbox}{max width=\linewidth}
$\displaystyle b=\frac5{5+4e},\quad \alpha=\frac5{4e-15},\quad
V=(1-z)^5,\quad N=(1-z/b)^{5-e}(1-z/\alpha)^e,\quad W=N/V,$
\end{adjustbox}
\end{center}
\begin{center}
\begin{adjustbox}{max width=\linewidth}
$\displaystyle A=(1-z/\alpha)(N-hV)/(4(1-h)),\quad
R=\frac{aW+5/e}{kW+1},\quad k=\frac{5-e}{5h},\quad a=-(4+5/e)k.$
\end{adjustbox}
\end{center}
The coefficient $p'(0)=c-1$ determines $h=(25\beta^2+K_e)/(25(\beta^2+1))$, with $K_e=36,39,34,21$. For real $\beta\ne0$, $h\ne0,1$. The four values of $W_{zz}(0)$ are $-64/5,-96/5,-96/5,-64/5$, so the required local ramification is nondegenerate. Coefficient comparison gives exactly the four original-parameter maps in the result.

\paragraph*{Boundaries, intersections and the final completeness statement}
The source list has 69 records. Eight root-label duplicates and ten whole-family inclusions reduce it to 51. Nine independent infinity boundaries are merged by $k_0=1/(h-1)$ and, where required, $u_0=t_0/(h-1)$; the boundary is $k_0=0$. This gives the 42 covering families in Theorem~\ref{thm:main}. The three special quadratic-factor intersections $d=-3,-1,1/9$ lie in degree-two conditions $\CCfam{9},\CCfam{10},\CCfam{11}$, respectively; an intersection does not justify deleting the containing type-II family.

Every apparent denominator has been checked in original parameters. Thirty-four source $h=1$ pole layers admit no finite coefficient vector: 32 have a nonzero constant Laurent residue and two have exact Bézout unit certificates. For example D3-I-01 has the unavoidable residue $420$, while D3-I-02 has $63$. The three zero-$W''(0)$ shapes $d=13,3,11/3$ have nonremovable original-coefficient residues $96r$, $12r$, $128r/9$, respectively, where $r=\beta^2\ne0$. Ten other auxiliary-chart denominator boundaries restore directly into F. Additional regular $k_0=0$ points also satisfy all five F equations. None is discarded merely because an auxiliary formula has a pole.

\paragraph*{Exact-computation requirements for independent reproduction}
The finite branch-elimination claims above are algebraic statements, not numerical tests. To make their use in the completeness argument independently auditable, the supplementary computational material should include: (i) the full integer partition and root-label enumeration with its exclusion rules; (ii) the 32 explicit nonzero Laurent residues and the two B\'ezout identities for the $h=1$ layers; (iii) all zero-$W''(0)$ and denominator-boundary substitutions in the ten original coefficients; and (iv) the coefficient-by-coefficient verification of the 27 parameter transports recorded in \ref{app:records}. The present TeX source lists the conclusions and their destinations but does not itself contain all of these executable certificates. The counting identity below is therefore a check of the claimed finite enumeration, not an independent substitute for these exact algebraic outputs.

For example, the source D3-I-01 is exactly D3-II-01 at $d=-35$, D3-I-09 is the same type-II family at $d=1$, and D3-I-10 is D3-II-02 at $d=0$; the remaining free $k$ is unchanged in each map. These are whole-family containments, not tests at chosen numerical parameters. All 27 duplicate, inclusion and boundary transports are identities of the ten original coefficients. The finite rational membership test, the degree bound, the exhaustive zero-parameter analysis, the full integer partitions and the restored boundaries together show that every real center satisfies one of $\CCfam{1},\ldots,\CCfam{42}$. The converse is proved independently below.

\paragraph*{Why the final count proves coverage and does not replace it}
There are four logically separate obligations. First, the analytic equivalence sends any center into the proportional layer or a full-field degree $2,3,4,5$. Second, the zero/discriminant cases and the finite integer partitions exhaust those possibilities. Third, solving the complete residue systems and restoring every coefficient gives the 69 source records. Fourth, the 27 exact transports and all denominator checks replace that source cover by the stated 42-family union. The count $69-8-10-9=42$ verifies the bookkeeping only after these four implications have been established. \ref{app:records} records the destinations of all source branches.

\medskip
\noindent\textbf{Source-to-family index.}
Table~\ref{tab:atlas} summarizes the source branches of the forty-two conditions in Theorem~\ref{thm:main}. The preceding classification and the transports in \ref{app:records} account for the omitted source records and their retained destinations.
\begin{longtable}{@{}>{\raggedright\arraybackslash}p{0.19\linewidth}>{\raggedright\arraybackslash}p{0.35\linewidth}>{\raggedright\arraybackslash}p{0.37\linewidth}@{}}
\caption{Source records and retained parameter families.}\label{tab:atlas}\\
\toprule Family & Exact source branch label & Sufficiency certificate class\\\midrule\endfirsthead
\toprule Family & Exact source branch label & Sufficiency certificate class\\\midrule\endhead
$\CCfam1$ & F & Integrating factor and first integrals\\
$\CCfam2$ & A & Original reflection\\
$\CCfam3$--$\CCfam4$ & G, H & Algebraic reversible realization\\
$\CCfam5$--$\CCfam7$ & D2-swap\_1, D2-swap\_2, D2-swap\_3 & Original-time nonlinear reverser\\
$\CCfam8$--$\CCfam{13}$ & D2-fixed\_1, \_3, \_4, \_5, \_6, \_7 & Original-field orbital reverser\\
$\CCfam{14}$--$\CCfam{22}$ & D3-I-02--08, D3-I-11--12 & Explicit $W_i,\ell_i,j_i,\Omega_i$\\
$\CCfam{23}$--$\CCfam{24}$ & D3-II-01--02 & Explicit $W_i,\ell_i,j_i,\Omega_i$\\
$\CCfam{25}$--$\CCfam{29}$ & D3-III-01,02,03,07,08 & Explicit $W_i,\ell_i,j_i,\Omega_i$\\
$\CCfam{30}$--$\CCfam{35}$ & D4-I-01,02,05,06,07,09 & Explicit $W_i,\ell_i,j_i,\Omega_i$\\
$\CCfam{36}$ & D4-II-01 & Mixed-fiber algebraic reversibility\\
$\CCfam{37}$--$\CCfam{38}$ & D4-III-01--02 & Two-positive-root fiber certificates\\
$\CCfam{39}$--$\CCfam{42}$ & D5-01--04 & Degree-five mixed-fiber certificates\\
\bottomrule
\end{longtable}

\section{Sufficiency and explicit reversible realizations}\label{sec:sufficiency}
We prove the conditions in their stated numerical order. Throughout this section $\delta=0$, $X_i=(P_i,Q_i)$ denotes \eqref{eq:family} with the coefficients of condition $\CCfam{i}$, and $P_{i,5},Q_{i,5}$ denote its homogeneous quintic parts. Set $g=a_6$, $c_i=a_1+1$, $s=y/x$, $p_i(s)=P_{i,5}(1,s)$, $A_i(s)=Q_{i,5}(1,s)-sp_i(s)$ and $B_i(s)=-1/4-gs+c_is^2$, consistently with \eqref{eq:AB}. Unsubscripted letters refer to the family currently under discussion. Parameter domains are exactly those stated in \ref{app:atlas}; no new free parameter is introduced in the proof.

\subsection{Conditions 1--4: integrating factors and elementary reversibility}
\paragraph*{Condition 1: explicit integrating factor and first integrals}\label{suff:1}
For \hyperref[center:1]{$\CCfam{1}$}, put $M=4a_1+2a_7$, $c=a_1+1$, $g=a_6$, $s=y/x$. The polynomials $p,A$ are those of the original system in \eqref{eq:AB}. The five equations $R_2=\cdots=R_6=0$ in condition $\CCfam{1}$ give \eqref{eq:Fidentity} without division by a parameter. Set $H_B=1+4gs-4cs^2=-4B$ and define
\begin{equation*}
\phi(s)=H_B(s)^{M/2}\exp\!\left(\frac{g(M+2)}2J_B(s)\right),\qquad
\rho(x,y)=x^M\phi(y/x).
\end{equation*}
Put $\Delta=g^2+c$. The following expression defines $J_B$ by $J_B(0)=0$ and $J_B'=1/B$:
{\footnotesize\begin{equation*}
J_B(s)=\begin{cases}
-4s,&c=g=0,\\
-g^{-1}\log(1+4gs),&c=0,\ g\ne0,\\
\displaystyle\frac1{\sqrt\Delta}\log\frac{(2cs-g-\sqrt\Delta)(-g+\sqrt\Delta)}{(2cs-g+\sqrt\Delta)(-g-\sqrt\Delta)},&c\ne0,\ \Delta>0,\\
\displaystyle-\frac2{2cs-g}-\frac2g,&c\ne0,\ \Delta=0,\\
\displaystyle\frac2{\sqrt{-\Delta}}\left(\arctan\frac{2cs-g}{\sqrt{-\Delta}}-\arctan\frac{-g}{\sqrt{-\Delta}}\right),&c\ne0,\ \Delta<0,
\end{cases}.
\end{equation*}}
All powers and logarithms are taken on the real analytic branches at $x=H_B=1$. In particular, $\rho(1,0)=1$. Its logarithmic derivative is $\phi'/\phi=(Mcs+g)/B$. The divergence calculation separates by homogeneity:
\[
\frac{\operatorname{div}(\rho X)}{x^M\phi}
=-Mcs-g+B\phi'/\phi
+x^4\bigl((M+6)p+A'+A\phi'/\phi\bigr)=0,
\]
where the last equality is \eqref{eq:Fidentity}. Thus $\rho$ is an explicit nonzero integrating factor in the original coordinates.

Here the first integral can be made more explicit than a general path quadrature. Put $A_1=x^2B(y/x)=xQ_{\rm lin}-yP_{\rm lin}$ and $A_5=x^6A(y/x)=xQ_5-yP_5$. For $M\ne-2,-6$ one has the elementary expression
\begin{equation*}
\mathcal H_1=\frac{\rho A_1+1/4}{M+2}+\frac{\rho A_5-1/4}{M+6}.
\end{equation*}
Euler's identity for the two closed homogeneous pieces of $\rho Q\,dx-\rho P\,dy$ proves $d\mathcal H_1=\rho Q\,dx-\rho P\,dy$. Both exceptional exponents can be treated directly, without an unproved limiting path in parameter space.

For $M=-2$, $\phi=1/H_B$ and a first integral is
\begin{equation}\label{eq:HFminus2}
\mathcal H_1=\frac{x^4A(s)}{4H_B(s)}-\frac1{16}
-\frac14\log x-\frac18\log H_B(s)-\frac g8J_B(s).
\end{equation}
For $M=-6$, identity \eqref{eq:Fidentity} gives $A'/A=(6cs-g)/B$, so $\phi=1/(4A)$. A first integral is
\begin{equation}\label{eq:HFminus6}
\mathcal H_1=-\frac{x^{-4}B(s)}{16A(s)}-\frac1{16}+\frac14\log x-\frac14\mathcal R_A(s).
\end{equation}
The last function also has a finite elementary expression. Let $\lambda$ run over the distinct complex roots of $A$, with multiplicities $m_\lambda$. Divide $p/A$ into its polynomial part $\sum_{j=0}^rq_js^j$ and partial fractions $\sum_{\lambda}\sum_{j=1}^{m_\lambda}b_{\lambda j}(s-\lambda)^{-j}$. All roots $\lambda$ of $A$ are nonzero. Then
\begin{equation*}
\begin{aligned}
\mathcal R_A(s)={}&\sum_{j=0}^r\frac{q_js^{j+1}}{j+1}
 +\sum_\lambda b_{\lambda1}\log(1-s/\lambda)\\
&+\sum_\lambda\sum_{j=2}^{m_\lambda}\frac{b_{\lambda j}}{1-j}
 \bigl((s-\lambda)^{1-j}-(-\lambda)^{1-j}\bigr),\\
b_{\lambda j}={}&\frac1{(m_\lambda-j)!}
\left[\frac{d^{m_\lambda-j}}{ds^{m_\lambda-j}}
 \left((s-\lambda)^{m_\lambda}\frac{p(s)}{A(s)}\right)\right]_{s=\lambda}.
\end{aligned}
\end{equation*}
The polynomial coefficients $q_j$ are obtained by ordinary Euclidean division. Empty sums are zero. Thus no unknown function remains: $\mathcal R_A(0)=0$, $\mathcal R_A'=p/A$. Conjugate roots and normalized logarithms give a real analytic result. This expression covers repeated roots and all degree drops of $A$, including constant $A$.

Direct differentiation verifies \eqref{eq:HFminus2} and \eqref{eq:HFminus6} against the same closed one-form. In every case $\mathcal H_1(1,0)=0$, $d\mathcal H_1(1,0)=0$ and $D^2\mathcal H_1(1,0)=I_2$. Small positive levels are closed curves, proving sufficiency for the entire F family.

\paragraph*{Condition 2: original reversibility}\label{suff:2}
Under \hyperref[center:2]{$\CCfam{2}$}, namely $a_2=a_4=a_6=a_8=a_{10}=0$, inspection of the normalized system shows that the first right-hand side is odd in $y$ and the second even. Hence $\mathcal J(x,y)=(x,-y)$ satisfies $X(\mathcal Jp)=-D\mathcal JX(p)$. Its fixed axis contains $(1,0)$ and the linear part has eigenvalues $\pm\ii$, proving a center without a coordinate or time change.

\paragraph*{Condition 3: an explicit radical change to a reversible field}\label{suff:3}
Here $T=a_2$, $U=a_3$ and $m=a_7-2$ are the parameters of \hyperref[center:3]{$\CCfam{3}$}; $U$ in this paragraph is a coefficient, not a phase coordinate. Use
\[
\eta=y\sqrt{1-2Ty/(3x)},\qquad
\frac{d\tau}{dt}=\frac{1-Ty/x}{\sqrt{1-2Ty/(3x)}}.
\]
Both the coordinate Jacobian and time factor equal one at $(1,0)$. Direct chain-rule substitution gives
\[
x_\tau=\eta(-x^4+Ux^2\eta^2),\qquad
\eta_\tau=(x^5-x)/4+(m+2)x^3\eta^2-\tfrac U2x\eta^4.
\]
This is reversible under $(x,\eta)\mapsto(x,-\eta)$ and has an elliptic linear part. The positive branch gives a center of the original system, including $T=0$ and $U=0$.

\paragraph*{Condition 4: remove a positive orbital multiplier}\label{suff:4}
The coefficients in \hyperref[center:4]{$\CCfam{4}$} imply the useful factorization
\[
\binom{\dot x}{\dot y}=(x+4gy)
\binom{-x^3y+a_3xy^3}{-1/4+x^4/4+a_7x^2y^2+a_9y^4}.
\]
With $d\tau/dt=x+4gy$, the remaining field is odd/even in $y$ and hence reversible under $(x,y)\mapsto(x,-y)$. The multiplier is one at the equilibrium and positive nearby. This proves sufficiency for every H parameter, including $g=0$.

\subsection{Conditions 5--13: nonlinear reversibility}
\label{sec:suff-low}
\paragraph*{Direct nonlinear symmetries of the original polynomial field}
Conditions 5--13 share a nonlinear reversal construction; for conditions 5--7 it preserves the original time.
For the original field write only for differentiation $s=y/x$, so that
\[\dot x/x=-cs+x^4p(s),\qquad \dot s=B(s)+x^4A(s).\]
Let $b$ be the constant specified in the corresponding certificate, set $\sigma=-s/(1+bs)$, and suppose the explicit rational functions $F,\rho$, with $F(0)=\rho(0)=1$, satisfy
\begin{center}
\begin{adjustbox}{max width=\linewidth}
$\displaystyle\begin{aligned}
&\rho=-\sigma'B/(B\circ\sigma), F=A(B\circ\sigma)/[B(A\circ\sigma)],\qquad -cs+\frac{F'}{4F}B=\rho c\sigma,p+\frac{F'}{4F}A=-\rho F(p\circ\sigma).
\end{aligned}$
\end{adjustbox}
\end{center}
Then $\mathcal R(x,y)=(xr(s),xr(s)\sigma(s))$, $r=F^{1/4}$, satisfies
\[\mathcal R^2=\mathrm{id},\qquad D\mathcal R\,X=-\rho(s)X\circ\mathcal R.\]
Indeed $F(F\circ\sigma)=1$ and $\sigma\circ\sigma=s$, giving the first identity. Differentiating $\log(xr)$ and $\sigma(s)$ along the displayed original field proves the two components of the second identity separately, by comparing their constant and $x^4$ coefficients. This is an identity of the original vector field, not a center criterion borrowed from a scalar equation. The nine rows displayed below supply explicit rational data for these identities, with all free parameters retained.

If $\rho=1$, this is a reverser in the original time. If $\rho\ne1$, $\rho(\sigma)=1/\rho$ and the positive time factor $a=(1+\rho)/2$ gives $a\circ\mathcal R=a/\rho$. Hence $Y=X/a$ is reversed by $\mathcal R$. In either case the explicit coordinates
\[\xi=(x+\mathcal R_1)/2,\qquad\eta=(y-\mathcal R_2)/2\]
satisfy $\xi\circ\mathcal R=\xi$, $\eta\circ\mathcal R=-\eta$. Their Jacobian is one at $(1,0)$: $D\mathcal R(1,0)=\left(\begin{smallmatrix}1&r'(0)\\0&-1\end{smallmatrix}\right)$. Thus this gives a local analytic reflection and proves the center by reversibility of an elliptic equilibrium. All statements are local near $(1,0)$; denominators and chosen roots there are regular.

For $i=5,\ldots,13$, each row specifies $b_i$, $F_i$ and $\rho_i$.
The same definitions apply to every row:
\begin{equation}\label{eq:app-low-common}
\begin{aligned}
\sigma_i(s)&=-\frac{s}{1+b_i s},\qquad
r_i(s)=F_i(s)^{1/4},\\
\mathcal R_i(x,y)&=\bigl(xr_i(y/x),xr_i(y/x)\sigma_i(y/x)\bigr).
\end{aligned}
\end{equation}
The branch of $F_i^{1/4}$ is the positive real-analytic branch at $s=0$.
Inserting the coefficients from \ref{app:atlas} in the original vector field $X_i$
gives the identities
\begin{equation}\label{eq:app-low-test}
F_i(s)F_i(\sigma_i(s))=1,\quad
\rho_i(\sigma_i(s))\rho_i(s)=1,\quad
D\mathcal R_i\,X_i=-\rho_i(s)X_i\circ\mathcal R_i.
\end{equation}
For conditions 5--7, $\rho_i=1$. For conditions 8--13, the positive factor $a_i=(1+\rho_i)/2$ gives the rescaled reversible field proved above.
The nine explicit rational triples $(b_i,F_i,\rho_i)$ for $5\le i\le13$ are recorded in \ref{app:suffdata}; they satisfy \eqref{eq:app-low-test} on their stated parameter domains.

\subsection{Conditions 14--42: algebraic changes to reversible systems}\label{sec:suff-high}
\paragraph*{A direct algebraic symmetrization lemma}
For the current family define $J_i,L_i$ by the following formulas, with subscripts suppressed:
\[
J(s)=\frac{cs}{B(s)}+\frac{p(s)}{A(s)}=\frac{D}{AB},\qquad
L(s)=\frac{B'}B-\frac{A'}A-\frac{4p}A.
\]
Choose an explicitly specified $\chi$ with $\chi(0)=0$, $\chi'(0)=1$, and set
\begin{equation}\label{eq:directchange}
U=x[-A(y/x)/B(y/x)]^{1/4},\qquad V=\chi(y/x),\qquad
\frac{d\tau}{dt}=-B(s)\chi'(s).
\end{equation}
All roots use the branch equal to one at $s=0$. The coordinate Jacobian at $(1,0)$ is one and the time factor is $1/4>0$. Direct logarithmic differentiation gives
\begin{equation}\label{eq:directsystem}
U_\tau=\frac{U}{\chi'(s)}\left[J(s)-\frac{U^4-1}{4}L(s)\right],\qquad
V_\tau=U^4-1.
\end{equation}
Indeed, $U^4=-x^4A/B$, $\dot s=B(1-U^4)$ and
\[
\frac{\dot U}{U}=\frac{\dot x}{x}
+\frac14\left(\frac{A'}A-\frac{B'}B\right)\dot s
=-B\left[J-\frac{U^4-1}{4}L\right].
\]
If $J/\chi'$ and $L/\chi'$ are odd functions of $V$, then \eqref{eq:directsystem} is reversible under $(U,V)\mapsto(U,-V)$. Its linearization is $\bigl(\begin{smallmatrix}0&-4\\4&0\end{smallmatrix}\bigr)$, because $J'(0)=-4$. Consequently it is a center, and the inverse coordinate and positive time changes give a center of the original field.

For $i=14,\ldots,42$, certificate $i$ below specifies the quantities
$c_i$, $\alpha_i$, $\Omega_i$, $\lambda_i$, $\kappa_i$, $\nu_i$, $\gamma_i$ and the
three polynomials $N_i(z),V_i(z),E_i(z)$. Their argument $z$ and the argument $w$ below are dummy scalar variables. The symbols $g,k,u,d$ retain exactly the values and domains in condition $\CCfam{i}$; $c_i=a_1+1$, $\alpha_i$ rescales $s$, and $\Omega_i=W_i^{\prime\prime}(0)$. The constants $\lambda_i,\kappa_i,\nu_i,\gamma_i$ are the displayed functions of those parameters, rather than additional choices. The denominator polynomial $V_i(z)$ is distinct from the new phase coordinate $V$.
\emph{All} the remaining row functions are defined, uniformly and without
unlisted choices, by
\begin{equation}\label{eq:app-high-defs}
\begin{aligned}
&B_i(s)=-\tfrac14-gs+c_i s^2,\qquad W_i(s)=\frac{N_i(\alpha_i s)}{V_i(\alpha_i s)},\\
&n_i(w)=\nu_i(w-1)+\gamma_i,\qquad j_i(w)=\frac{(\lambda_i+1)n_i(w)}{4w(\kappa_iw-\kappa_i-1)},\\
&\ell_i(w)=-\frac{(\lambda_i+1)\bigl[n_i(w)+\kappa_i w-
\lambda_i(\kappa_i+1)/(\lambda_i+1)\bigr]}{w(\kappa_iw-\kappa_i-1)},\\
&A_i(s)=\frac{E_i(\alpha_i s)}4
\bigl[(\kappa_i+1)V_i(\alpha_i s)-\kappa_iN_i(\alpha_i s)\bigr],\\
&p_i(s)=A_i(s)\left[j_i(W_i(s))W_i'(s)-\frac{c_i s}{B_i(s)}\right].
\end{aligned}
\end{equation}
The $\lambda_i+1$ factors are nonzero for every listed row. The denominators
in these formulas are regular at $s=0$ (in particular
$W_i(0)=1$ and $\kappa_i W_i(0)-\kappa_i-1=-1$).

For conditions 14--42 the row data in \ref{app:suffdata} give explicit rational functions satisfying
\begin{equation}\label{eq:highpullback}
\begin{aligned}
&W_i(0)=1,\quad W_i'(0)=0,\quad \Omega_i=W_i''(0)\ne0,\\
&L_i=\ell_i(W_i)W_i',\quad J_i=j_i(W_i)W_i'.
\end{aligned}
\end{equation}
Take
\begin{equation}\label{eq:highchi}
\chi_i(s)=s\sqrt{\frac{2(W_i(s)-1)}{\Omega_i s^2}}.
\end{equation}
Then $W_i=1+\Omega_iV^2/2$ and $W_i'=\Omega_iV\chi_i'$. The actual transformed system is
\begin{equation}\label{eq:highsystem}
\begin{aligned}
U_\tau&=\Omega_iUV\left[j_i\!\left(1+\frac{\Omega_iV^2}{2}\right)
 -\frac{U^4-1}{4}\ell_i\!\left(1+\frac{\Omega_iV^2}{2}\right)\right],\\
V_\tau&=U^4-1.
\end{aligned}
\end{equation}
Its first component is visibly odd in $V$ and the second even. No inverse of $W_i$ and no solution of an involution equation is required. The certificates in \ref{app:suffdata} supply the row-dependent rational data; the original vector field is specified once by the corresponding condition in \ref{app:atlas}.

Let $A_{ij}$ and $p_{ij}$ denote the coefficients of $s^j$ in the original polynomials $A_i$ and $p_i$, respectively. For every row the required identities are
\begin{equation}\label{eq:app-high-test}
\begin{gathered}
A_i(s)=\sum_{j=0}^{6}A_{ij}s^j,\qquad
p_i(s)=\sum_{j=0}^{5}p_{ij}s^j,\\
\text{the coefficients agree with the corresponding family in \ref{app:atlas}},\\
W_i(0)=1,\quad W_i'(0)=0,\quad W_i''(0)=\Omega_i\ne0,\\
L_i(s)=\ell_i(W_i(s))W_i'(s),\qquad
J_i(s)=j_i(W_i(s))W_i'(s).
\end{gathered}
\end{equation}
Thus \eqref{eq:app-high-test} identifies the reconstruction with the original coefficients and verifies the pullbacks used in \eqref{eq:highsystem}. All identities are understood on the full stated parameter domains, including their retained boundary strata.
The twenty-nine explicit rational-data rows required in \eqref{eq:app-high-defs}--\eqref{eq:app-high-test} are recorded in \ref{app:suffdata}. The proof of the unified criterion is the calculation above, so the same change of variables is not repeated for every row.

\subsection{Completion of the sufficiency proof}
\begin{theorem}[Sufficiency of the reconstructed conditions]\label{thm:suff}
Each of the forty-two original-coefficient conditions in Theorem~\ref{thm:main}
is sufficient for an elementary center at $(1,0)$.
\end{theorem}
\begin{proof}
For $\CCfam1$ use the integrating factor and first integrals constructed
above, including the exceptional exponents $M=-2,-6$. For $\CCfam2$
use the original reflection; for $\CCfam3$--$\CCfam4$ use the explicit
positive-time analytic transformations to reversible systems. For
$\CCfam5$--$\CCfam{13}$, apply the original or orbital reversal identities
\eqref{eq:app-low-test}. For $\CCfam{14}$--$\CCfam{42}$, apply the
uniform pullback and reversal argument \eqref{eq:highpullback}--
\eqref{eq:highsystem}, with the row data in \ref{app:suffdata} and coefficient checks
\eqref{eq:app-high-defs}--\eqref{eq:app-high-test}. For conditions 2--42 the resulting reflection fixes an elementary elliptic equilibrium, so nearby trajectories are closed; for condition 1 the same conclusion follows from the positive-definite Hessian of $\mathcal H_1$. Inverting the local coordinate changes and
positive time factors preserves the center property.
\end{proof}

The necessity part of Theorem~\ref{thm:main} is established in
Section~\ref{sec:classification}, using the exact coefficient transports
in \ref{app:records}; the sufficiency part is Theorem~\ref{thm:suff}.

\appendix
\section{The complete original-parameter atlas}\label{app:atlas}
All parameters not explicitly determined are free real parameters subject to the domains shown in each row; $g=a_6$. Labels $\CCfam1$--$\CCfam{42}$ designate covering families, and a boundary point can belong to more than one of them. Auxiliary numerator polynomials are defined where used and are not additional free variables.

\begingroup
\normalfont
\scriptsize
\allowdisplaybreaks[4]
\setlength{\jot}{1pt}
\setlength{\abovedisplayskip}{4pt}
\setlength{\belowdisplayskip}{4pt}
\phantomsection\paragraph*{Condition $\CCfam{1}$ (F)}\label{center:1}
Set $c=a_1+1$, $g=a_6$, $M=4a_1+2a_7$, and define the finite coefficient lists
\begin{align*}
&(A_0,\ldots,A_6)=(1/4,g,a_7-a_1,a_8-a_2,a_9-a_3,a_{10}-a_4,-a_5),\\
&\qquad (p_0,\ldots,p_5)=(0,a_1,a_2,a_3,a_4,a_5).
\end{align*}
All unlisted subscripts are zero. The five short original-parameter conditions are $R_k=0$ for $k=2,\ldots,6$, where
\begin{align*}
&R_k={}-\frac{k+1}{4}A_{k+1}+g(1-k)A_k+c(M+k-1)A_{k-1}+(M+6)(-p_k/4-gp_{k-1}+cp_{k-2}).
\end{align*}
\phantomsection\paragraph*{Condition $\CCfam{2}$ (A)}\label{center:2}
Free parameters: $a_{1}, a_{3}, a_{5}, a_{7}, a_{9}$.
\[a_2=a_4=a_6=a_8=a_{10}=0.\]
\phantomsection\paragraph*{Condition $\CCfam{3}$ (G)}\label{center:3}
Free parameters: $T, U, m$.
\begin{align*}
&a_{1}=-1,\quad a_{2}=T,\quad a_{3}=U,\quad a_{4}=- \frac{5 T U}{3},\quad a_{5}=\frac{2 T^{2} U}{3},\\
&a_{6}=0,\quad a_{7}=m + 2,\quad a_{8}=- \frac{T \left(2 m + 3\right)}{3},\quad a_{9}=- \frac{U}{2},\quad a_{10}=\frac{T U}{3}.
\end{align*}
\phantomsection\paragraph*{Condition $\CCfam{4}$ (H)}\label{center:4}
Free parameters: $g, a_{3}, a_{7}, a_{9}$.
\begin{align*}
&a_{1}=-1,\quad a_{2}=- 4 g,\quad a_{4}=4 a_{3} g,\quad a_{5}=0,\quad a_{6}=g,\quad a_{8}=4 a_{7} g,\quad a_{10}=4 a_{9} g.
\end{align*}
\phantomsection\paragraph*{Condition $\CCfam{5}$ (D2-swap\_1)}\label{center:5}
Free parameters: $g, a_{3}, a_{4}, a_{5}$.
Domain: $g\ne0$.
\begin{align*}
&a_{1}=15 g^{2} - 1,\quad a_{2}=2 g,\quad a_{6}=g,\quad a_{7}=\frac{6 a_{3} g + a_{4} + 2400 g^{5} - 264 g^{3}}{120 g^{3}},\\
&a_{8}=- 2 g \left(20 g^{2} - 1\right),\quad a_{9}=- \frac{24 a_{3} g^{2} + 14 a_{4} g - 5 a_{5} + 23040 g^{6} - 2016 g^{4}}{60 g^{2}},\\
&a_{10}=- \frac{36 a_{3} g^{2} - 9 a_{4} g + 5 a_{5} + 8640 g^{6} - 864 g^{4}}{15 g}.
\end{align*}
\phantomsection\paragraph*{Condition $\CCfam{6}$ (D2-swap\_2)}\label{center:6}
Free parameters: $g, a_{4}, a_{5}$.
Domain: $g\ne0$.
\begin{align*}
&a_{1}=3 g^{2} - 1,\quad a_{2}=2 g,\quad a_{3}=- \frac{2 a_{4} g + a_{5} + 48 g^{6}}{4 g^{2}},\quad a_{6}=g,\\
&a_{7}=- \frac{4 a_{4} g + 3 a_{5} - 48 g^{6} + 80 g^{4}}{48 g^{4}},\quad a_{8}=- 2 g \left(2 g - 1\right) \left(2 g + 1\right),\\
&a_{9}=\frac{2 a_{4} g + a_{5} - 32 g^{6} + 32 g^{4}}{4 g^{2}},\quad a_{10}=- \frac{a_{4} - 48 g^{5} + 32 g^{3}}{3}.
\end{align*}
\phantomsection\paragraph*{Condition $\CCfam{7}$ (D2-swap\_3)}\label{center:7}
Free parameters: $g, a_{5}$.
Domain: $g\ne0$.
\begin{align*}
&a_{1}=\frac{7 g^{2} - 9}{9},\quad a_{2}=2 g,\quad a_{3}=\frac{243 a_{5} - 112 g^{6} - 144 g^{4}}{108 g^{2}},\\
&a_{4}=- \frac{729 a_{5} - 112 g^{6} - 72 g^{4}}{243 g},\quad a_{6}=g,\quad a_{7}=\frac{243 a_{5} + 16 g^{6} - 144 g^{4}}{144 g^{4}},\\
&a_{8}=- \frac{2 g \left(28 g^{2} - 27\right)}{27},\quad a_{9}=- \frac{243 a_{5} - 128 g^{6} + 144 g^{4}}{108 g^{2}},\quad a_{10}=\frac{243 a_{5} - 32 g^{6} + 72 g^{4}}{243 g}.
\end{align*}
\phantomsection\paragraph*{Condition $\CCfam{8}$ (D2-fixed\_1)}\label{center:8}
Free parameters: $g, a_{4}, a_{5}$.
Domain: $g\ne0$.
\begin{align*}
&a_{1}=3 g^{2} - 1,\quad a_{2}=- 10 g,\quad a_{3}=\frac{6 a_{4} g - a_{5} - 3888 g^{6} - 864 g^{4}}{36 g^{2}},\quad a_{6}=g,\\
&a_{7}=\frac{24 a_{4} g - 7 a_{5} - 18576 g^{6} - 3888 g^{4}}{2160 g^{4}},\quad a_{8}=\frac{24 a_{4} g - 7 a_{5} - 7776 g^{6} - 4428 g^{4}}{270 g^{3}},\\
&a_{9}=\frac{2 a_{4} g - a_{5} + 2592 g^{6} - 864 g^{4}}{60 g^{2}},\quad a_{10}=- \frac{9 a_{4} g - 2 a_{5} - 1296 g^{6} - 1728 g^{4}}{15 g}.
\end{align*}
\phantomsection\paragraph*{Condition $\CCfam{9}$ (D2-fixed\_3)}\label{center:9}
Free parameters: $g, a_{3}, a_{4}, a_{5}$.
Domain: $g\ne0$.
\begin{align*}
&a_{1}=- \frac{5 g^{2} + 9}{9},\quad a_{2}=- 2 g,\quad a_{6}=g,\quad a_{7}=- \frac{162 a_{3} g - 243 a_{4} + 800 g^{5} + 1224 g^{3}}{1080 g^{3}},\\
&a_{8}=- \frac{162 a_{3} g - 243 a_{4} + 1400 g^{5} + 954 g^{3}}{405 g^{2}},\quad a_{9}=\frac{8424 a_{3} g^{2} - 5346 a_{4} g + 729 a_{5} - 17920 g^{6} + 3168 g^{4}}{4860 g^{2}},\\
&a_{10}=\frac{1620 a_{3} g^{2} - 1215 a_{4} g + 243 a_{5} - 1600 g^{6} + 1440 g^{4}}{1215 g}.
\end{align*}
\phantomsection\paragraph*{Condition $\CCfam{10}$ (D2-fixed\_4)}\label{center:10}
Free parameters: $g, a_{3}, a_{4}, a_{5}$.
Domain: $g\ne0$.
\begin{align*}
&a_{1}=- \frac{3 g^{2} + 4}{4},\quad a_{2}=- g,\quad a_{6}=g,\quad a_{7}=\frac{a_{3} g + a_{4} - 3 g^{5} - 5 g^{3}}{3 g^{3}},\\
&a_{8}=\frac{2 a_{3} g + 2 a_{4} - 9 g^{5} - 7 g^{3}}{3 g^{2}},\quad a_{9}=\frac{2 a_{3} g^{2} - a_{4} g + 3 a_{5} + 2 g^{4}}{3 g^{2}},\quad a_{10}=- \frac{2 a_{3} g - a_{4} - 6 g^{5} - 4 g^{3}}{3}.
\end{align*}
\phantomsection\paragraph*{Condition $\CCfam{11}$ (D2-fixed\_5)}\label{center:11}
Free parameters: $g, a_{4}, a_{5}$.
Domain: $g\ne0$.
\begin{align*}
&a_{1}=- \frac{21 g^{2} + 25}{25},\quad a_{2}=- \frac{2 g}{5},\quad a_{3}=- \frac{6250 a_{4} g + 15625 a_{5} - 336 g^{6} - 800 g^{4}}{2500 g^{2}},\quad a_{6}=g,\\
&a_{7}=- \frac{25000 a_{4} g + 78125 a_{5} + 7728 g^{6} + 13200 g^{4}}{8400 g^{4}},\quad a_{8}=- \frac{5000 a_{4} g + 15625 a_{5} + 2016 g^{6} + 1380 g^{4}}{1050 g^{3}},\\
&a_{9}=\frac{143750 a_{4} g + 203125 a_{5} + 74592 g^{6} + 69600 g^{4}}{52500 g^{2}},\quad a_{10}=- \frac{625 a_{4} g - 6250 a_{5} + 1008 g^{6} + 960 g^{4}}{4375 g}.
\end{align*}
\phantomsection\paragraph*{Condition $\CCfam{12}$ (D2-fixed\_6)}\label{center:12}
Free parameters: $g, a_{8}$.
Domain: $g\ne0$.
\begin{align*}
&a_{1}=- \frac{8 g^{2} + 9}{9},\quad a_{2}=0,\quad a_{3}=0,\quad a_{4}=0,\quad a_{5}=0,\\
&a_{6}=g,\quad a_{7}=\frac{3 a_{8} - 4 g}{4 g},\quad a_{9}=0,\quad a_{10}=0.
\end{align*}
\phantomsection\paragraph*{Condition $\CCfam{13}$ (D2-fixed\_7)}\label{center:13}
Free parameters: $g$.
Domain: $g\ne0$.
\begin{align*}
&a_{1}=- \frac{45 g^{2} + 49}{49},\quad a_{2}=\frac{2 g}{7},\quad a_{3}=\frac{36 g^{2} \left(30 g^{2} + 49\right)}{2401},\quad a_{4}=\frac{8 g^{3} \left(360 g^{2} + 539\right)}{16807},\\
&a_{5}=\frac{16 g^{4} \left(135 g^{2} + 196\right)}{117649},\quad a_{6}=g,\quad a_{7}=- \frac{134 g^{2} + 245}{49},\\
&a_{8}=- \frac{2 g \left(68 g^{2} + 105\right)}{49},\quad a_{9}=- \frac{4 g^{2} \left(499 g^{2} + 735\right)}{2401},\quad a_{10}=- \frac{40 g^{3} \left(34 g^{2} + 49\right)}{16807}.
\end{align*}
\phantomsection\paragraph*{Condition $\CCfam{14}$ (D3-I-02)}\label{center:14}
Free parameters: $g, k, u$.
Domain: $g\ne0$.
\begin{align*}
&a_{1}=- \frac{3 g^{2} + 4}{4},\quad a_{2}=- 3 g,\quad a_{3}=- \frac{3 g^{2} \left(6 g^{2} k + 24 g^{2} u - 3 g^{2} + 2\right)}{2},\\
&a_{4}=- g^{3} \left(- 6 g^{2} k + 84 g^{2} u - 6 g^{2} + 1\right),\quad a_{5}=\frac{9 g^{6} \left(28 k + 1\right)}{4},\\
&a_{6}=g,\quad a_{7}=\frac{6 g^{2} k - 3 g^{2} - 2}{2},\quad a_{8}=g \left(10 g^{2} k - 8 g^{2} - 3\right),\\
&a_{9}=- \frac{g^{2} \left(116 g^{2} k + 144 g^{2} u + 35 g^{2} + 12\right)}{4},\quad a_{10}=- g^{3} \left(84 g^{2} k + 84 g^{2} u + 3 g^{2} + 1\right).
\end{align*}
\phantomsection\paragraph*{Condition $\CCfam{15}$ (D3-I-03)}\label{center:15}
Free parameters: $g, k$.
Domain: $g\ne0$.
\begin{align*}
&a_{1}=- \frac{21 g^{2} + 25}{25},\quad a_{2}=- 2 g,\quad a_{3}=\frac{4 g^{2} \left(720 g^{2} k + 1854 g^{2} + 725\right)}{625},\\
&a_{4}=\frac{24 g^{3} \left(416 g^{2} k + 360 g^{2} + 175\right)}{625},\quad a_{5}=\frac{144 g^{6} \left(272 k + 27\right)}{3125},\\
&a_{6}=g,\quad a_{7}=\frac{48 g^{2} k - 6 g^{2} - 25}{25},\quad a_{8}=\frac{2 g \left(64 g^{2} k - 36 g^{2} - 25\right)}{25},\\
&a_{9}=- \frac{4 g^{2} \left(16 g^{2} k - 909 g^{2} - 725\right)}{625},\quad a_{10}=- \frac{24 g^{3} \left(544 g^{2} k - 1206 g^{2} - 875\right)}{3125}.
\end{align*}
\phantomsection\paragraph*{Condition $\CCfam{16}$ (D3-I-04)}\label{center:16}
Free parameters: $g, k$.
Domain: $g\ne0$.
\begin{align*}
&a_{1}=\frac{7 g^{2} - 9}{9},\quad a_{2}=- \frac{10 g}{3},\quad a_{3}=\frac{4 g^{2} \left(1680 g^{2} k - 266 g^{2} + 117\right)}{81},\\
&a_{4}=\frac{56 g^{3} \left(800 g^{2} k + 72 g^{2} - 9\right)}{243},\quad a_{5}=- \frac{3920 g^{6} \left(176 k + 1\right)}{729},\\
&a_{6}=g,\quad a_{7}=\frac{48 g^{2} k - 38 g^{2} - 9}{9},\quad a_{8}=\frac{2 g \left(256 g^{2} k - 156 g^{2} - 45\right)}{27},\\
&a_{9}=- \frac{4 g^{2} \left(912 g^{2} k - 397 g^{2} - 117\right)}{81},\quad a_{10}=- \frac{56 g^{3} \left(352 g^{2} k + 30 g^{2} + 9\right)}{243}.
\end{align*}
\phantomsection\paragraph*{Condition $\CCfam{17}$ (D3-I-05)}\label{center:17}
Free parameters: $g, k, u$.
Domain: $g\ne0$.
\begin{align*}
&a_{1}=- \frac{9 g^{2} + 25}{25},\quad a_{2}=- \frac{6 g}{5},\quad a_{3}=- \frac{12 g^{2} \left(144 g^{2} k + 576 g^{2} u - 18 g^{2} + 25\right)}{625},\\
&a_{4}=\frac{8 g^{3} \left(2016 g^{2} k + 2880 g^{2} u + 72 g^{2} - 25\right)}{3125},\quad a_{5}=- \frac{432 g^{6} \left(80 k - 1\right)}{15625},\\
&a_{6}=g,\quad a_{7}=\frac{48 g^{2} k - 6 g^{2} - 25}{25},\quad a_{8}=\frac{2 g \left(256 g^{2} k - 68 g^{2} - 75\right)}{125},\\
&a_{9}=- \frac{4 g^{2} \left(2000 g^{2} k + 1728 g^{2} u + 83 g^{2} + 75\right)}{625},\quad a_{10}=\frac{8 g^{3} \left(480 g^{2} k + 576 g^{2} u - 6 g^{2} - 5\right)}{625}.
\end{align*}
\phantomsection\paragraph*{Condition $\CCfam{18}$ (D3-I-06)}\label{center:18}
Free parameters: $g, k$.
Domain: $g\ne0$.
\begin{align*}
&a_{1}=- \frac{33 g^{2} + 49}{49},\quad a_{2}=- \frac{2 g}{7},\quad a_{3}=\frac{4 g^{2} \left(144 g^{2} k + 342 g^{2} + 245\right)}{2401},\\
&a_{4}=\frac{24 g^{3} \left(160 g^{2} k + 72 g^{2} - 49\right)}{16807},\quad a_{5}=\frac{144 g^{6} \left(16 k + 27\right)}{117649},\quad a_{6}=g,\\
&a_{7}=\frac{48 g^{2} k + 42 g^{2} - 49}{49},\quad a_{8}=\frac{2 g \left(256 g^{2} k + 180 g^{2} - 49\right)}{343},\\
&a_{9}=\frac{4 g^{2} \left(368 g^{2} k + 477 g^{2} + 245\right)}{2401},\quad a_{10}=\frac{24 g^{3} \left(32 g^{2} k + 18 g^{2} - 49\right)}{16807}.
\end{align*}
\phantomsection\paragraph*{Condition $\CCfam{19}$ (D3-I-07)}\label{center:19}
Free parameters: $g, k$.
Domain: $g\ne0$.
\begin{align*}
&a_{1}=- \frac{8 g^{2} + 9}{9},\quad a_{2}=- \frac{4 g}{3},\quad a_{3}=\frac{16 g^{2} \left(24 g^{2} k + 16 g^{2} + 9\right)}{81},\quad a_{4}=\frac{64 g^{3} \left(16 g^{2} k + 9\right)}{243},\\
&a_{5}=- \frac{2048 g^{6} \left(5 k + 1\right)}{729},\quad a_{6}=g,\quad a_{7}=\frac{12 g^{2} k - 8 g^{2} - 9}{9},\quad a_{8}=\frac{4 g \left(20 g^{2} k - 24 g^{2} - 9\right)}{27},\\
&a_{9}=- \frac{16 g^{2} \left(12 g^{2} k - 4 g^{2} - 9\right)}{81},\quad a_{10}=- \frac{64 g^{3} \left(20 g^{2} k - 12 g^{2} - 9\right)}{243}.
\end{align*}
\phantomsection\paragraph*{Condition $\CCfam{20}$ (D3-I-08)}\label{center:20}
Free parameters: $g, k$.
Domain: $g\ne0$.
\begin{align*}
&a_{1}=- \frac{45 g^{2} + 49}{49},\quad a_{2}=- \frac{6 g}{7},\quad a_{3}=\frac{180 g^{2} \left(16 g^{2} k + 70 g^{2} + 49\right)}{2401},\\
&a_{4}=\frac{24 g^{3} \left(1120 g^{2} k - 2120 g^{2} - 1323\right)}{16807},\quad a_{5}=- \frac{720 g^{6} \left(368 k + 125\right)}{117649},\\
&a_{6}=g,\quad a_{7}=\frac{48 g^{2} k - 6 g^{2} - 49}{49},\quad a_{8}=\frac{2 g \left(320 g^{2} k - 220 g^{2} - 147\right)}{343},\\
&a_{9}=- \frac{4 g^{2} \left(656 g^{2} k - 2125 g^{2} - 2205\right)}{2401},\quad a_{10}=- \frac{24 g^{3} \left(736 g^{2} k + 1870 g^{2} + 1323\right)}{16807}.
\end{align*}
\phantomsection\paragraph*{Condition $\CCfam{21}$ (D3-I-11)}\label{center:21}
Free parameters: $g$.
Domain: $g\ne0$.
\begin{align*}
&a_{1}=- \frac{77 g^{2} + 81}{81},\quad a_{2}=- \frac{2 g}{9},\quad a_{3}=\frac{44 g^{2} \left(938 g^{2} + 891\right)}{6561},\quad a_{4}=- \frac{8 g^{3} \left(6776 g^{2} + 6345\right)}{19683},\\
&a_{5}=\frac{16 g^{4} \left(10087 g^{2} + 9396\right)}{531441},\quad a_{6}=g,\quad a_{7}=\frac{1018 g^{2} + 891}{81},\\
&a_{8}=\frac{2 g \left(2228 g^{2} + 2133\right)}{243},\quad a_{9}=- \frac{4 g^{2} \left(16499 g^{2} + 15471\right)}{6561},\quad a_{10}=\frac{8 g^{3} \left(2926 g^{2} + 2727\right)}{19683}.
\end{align*}
\phantomsection\paragraph*{Condition $\CCfam{22}$ (D3-I-12)}\label{center:22}
Free parameters: $g$.
Domain: $g\ne0$.
\begin{align*}
&a_{1}=- \frac{65 g^{2} + 81}{81},\quad a_{2}=\frac{2 g}{9},\quad a_{3}=\frac{4 g^{2} \left(910 g^{2} + 1539\right)}{6561},\quad a_{4}=\frac{8 g^{3} \left(520 g^{2} + 837\right)}{19683},\\
&a_{5}=\frac{16 g^{4} \left(715 g^{2} + 1134\right)}{531441},\quad a_{6}=g,\quad a_{7}=- \frac{302 g^{2} + 567}{81},\\
&a_{8}=- \frac{2 g \left(508 g^{2} + 837\right)}{243},\quad a_{9}=- \frac{4 g^{2} \left(2075 g^{2} + 3321\right)}{6561},\quad a_{10}=- \frac{8 g^{3} \left(290 g^{2} + 459\right)}{19683}.
\end{align*}.
\phantomsection\paragraph*{Condition $\CCfam{23}$ (D3-II-01)}\label{center:23}
Free parameters: $d, g, k$.
Domain: $d - 13\ne0,\ g\ne0$.
\begin{align*}
&a_{1}=- \frac{5 g^{2} + 9}{9},\quad a_{2}=- \frac{2 g \left(3 d - 7\right)}{d - 13},\quad a_{3}=\frac{\mathcal N_{\mathrm{a},233}}{243 d - 3159},\\
&\mathcal N_{\mathrm{a},233}=- 2340 g^{4}+324 g^{2}- 11830 g^{4} k
       - 1620 d g^{2}\\
&- 60 d^{2} g^{4}+960 d g^{4}- 330 d^{2} g^{4} k+10 d^{3} g^{4} k+3510 d g^{4} k,\quad a_{4}=\frac{\mathcal N_{\mathrm{a},234}}{729 d - 9477},\\
&\mathcal N_{\mathrm{a},234}=- 3120 g^{5}- 1080 g^{3}- 1820 g^{5} k- 1512 d g^{3}\\
&- 240 d^{2} g^{5}+3360 d g^{5}- 5140 d^{2} g^{5} k+220 d^{3} g^{5} k+29780 d g^{5} k,\quad a_{5}=\frac{80 d g^{6} \left(5 d k - k - 1\right)}{729},\\
&a_{6}=g,\quad a_{7}=- \frac{d g^{2} k - d g^{2} - 13 g^{2} k + 3 g^{2} + 9}{9},\quad a_{8}=\frac{\mathcal N_{\mathrm{a},238}}{27 d - 351},\\
&\mathcal N_{\mathrm{a},238}=378 g+520 g^{3}- 1040 g^{3} k- 248 d g^{3}- 162 d g+16 d^{2} g^{3}- 16 d^{2} g^{3} k+288 d g^{3} k,\\
&a_{9}=\frac{\mathcal N_{\mathrm{a},239}}{243 d - 3159},\quad \mathcal N_{\mathrm{a},239}=324 g^{2}+1872 g^{4}- 2172 d g^{4}- 1620 d g^{2}\\
&+156 d^{2} g^{4}
       +7982 g^{4} k- 78 d^{2} g^{4} k- 2 d^{3} g^{4} k+738 d g^{4} k,\quad a_{10}=\frac{\mathcal N_{\mathrm{a},2310}}{729 d - 9477},\\
&\mathcal N_{\mathrm{a},2310}=- 1080 g^{3}- 1872 d g^{5}- 1512 d g^{3}+144 d^{2} g^{5}+1300 g^{5} k
       - 6652 d g^{5} k- 20 d^{3} g^{5} k+764 d^{2} g^{5} k.
\end{align*}
\phantomsection\paragraph*{Condition $\CCfam{24}$ (D3-II-02)}\label{center:24}
Free parameters: $d, g, k$.
Domain: $d - 3\ne0,\ g\ne0$.
\begin{align*}
&a_{1}=- \frac{3 g^{2} + 4}{4},\quad a_{2}=- \frac{4 d g}{d - 3},\quad a_{3}=\frac{d g^{2} \left(d g^{2} k + d g^{2} - 3 g^{2} k - 9 g^{2} - 4\right)}{4},\\
&a_{4}=\frac{\mathcal N_{\mathrm{a},244}}{d - 3},\quad \mathcal N_{\mathrm{a},244}=d^{3} g^{5}- 4 d^{2} g^{3}- 4 d^{2} g^{5}+3 d g^{5}+d^{3} g^{5} k- 4 d^{2} g^{5} k+3 d g^{5} k,\\
&a_{5}=\frac{d g^{6} \left(3 d k - k - 1\right)}{4},\quad a_{6}=g,\quad a_{7}=- \frac{d g^{2} k - d g^{2} - 3 g^{2} k - 3 g^{2} + 4}{4},\\
&a_{8}=\frac{\mathcal N_{\mathrm{a},248}}{d - 3},\quad \mathcal N_{\mathrm{a},248}=- g \left(d^{2} g^{2} k - d^{2} g^{2} - 4 d g^{2} k + 2 d g^{2} + 4 d + 3 g^{2} k + 3 g^{2}\right),\\
&a_{9}=\frac{g^{2} \left(d^{2} g^{2} - 3 d g^{2} k - 3 d g^{2} - 4 d + g^{2} k + g^{2}\right)}{4},\quad a_{10}=\frac{d^{2} g^{3} \left(d g^{2} - 3 g^{2} - 4\right)}{d - 3}.
\end{align*}
\phantomsection\paragraph*{Condition $\CCfam{25}$ (D3-III-01)}\label{center:25}
Free parameters: $g, k$.
Domain: $g\ne0$.
\begin{align*}
&a_{1}=\frac{11 g^{2} - 25}{25},\quad a_{2}=\frac{6 g}{5},\quad a_{3}=- \frac{12 g^{2} \left(48 g^{2} k + 22 g^{2} + 25\right)}{625},\quad a_{4}=\frac{8 g^{3} \left(96 g^{2} k + 88 g^{2} + 25\right)}{3125},\\
&a_{5}=\frac{48 g^{6} \left(16 k - 11\right)}{15625},\quad a_{6}=g,\quad a_{7}=- \frac{48 g^{2} k + 6 g^{2} + 25}{25},\quad a_{8}=- \frac{2 g \left(64 g^{2} k + 52 g^{2} - 75\right)}{125},\\
&a_{9}=\frac{4 g^{2} \left(80 g^{2} k + 77 g^{2} - 75\right)}{625},\quad a_{10}=\frac{8 g^{3} \left(32 g^{2} k - 30 g^{2} + 25\right)}{3125}.
\end{align*}
\phantomsection\paragraph*{Condition $\CCfam{26}$ (D3-III-02)}\label{center:26}
Free parameters: $g, k, u$.
Domain: $g\ne0$.
\begin{align*}
&a_{1}=8 g^{2} - 1,\quad a_{2}=0,\quad a_{3}=- 48 g^{2} \left(8 g^{2} k + 8 g^{2} u + 8 g^{2} - 1\right),\\
&a_{4}=128 g^{3} \left(8 g^{2} k + 8 g^{2} u + 8 g^{2} - 1\right),\quad a_{5}=0,\quad a_{6}=g,\quad a_{7}=- 12 g^{2} k - 12 g^{2} - 1,\\
&a_{8}=- 16 g^{3} \left(k + 1\right),\quad a_{9}=- 16 g^{2} \left(- 8 g^{2} k + 24 g^{2} u - 8 g^{2} - 3\right),\quad a_{10}=128 g^{3} \left(8 g^{2} u - 1\right).
\end{align*}
\phantomsection\paragraph*{Condition $\CCfam{27}$ (D3-III-03)}\label{center:27}
Free parameters: $g, k$.
Domain: $g\ne0$.
\begin{align*}
&a_{1}=35 g^{2} - 1,\quad a_{2}=6 g,\quad a_{3}=- 180 g^{2} \left(16 g^{2} k + 42 g^{2} - 1\right),\quad a_{4}=- 40 g^{3} \left(96 g^{2} k - 280 g^{2} + 35\right),\\
&a_{5}=1200 g^{6} \left(272 k + 245\right),\quad a_{6}=g,\quad a_{7}=- 48 g^{2} k - 54 g^{2} - 1,\quad a_{8}=2 g \left(64 g^{2} k - 68 g^{2} + 3\right),\\
&a_{9}=4 g^{2} \left(1424 g^{2} k + 605 g^{2} + 45\right),\quad a_{10}=- 40 g^{3} \left(544 g^{2} k - 210 g^{2} + 35\right).
\end{align*}
\phantomsection\paragraph*{Condition $\CCfam{28}$ (D3-III-07)}\label{center:28}
Free parameters: $g, k$.
Domain: $g\ne0$.
\begin{align*}
&a_{1}=\frac{5 g^{2} - 4}{4},\quad a_{2}=g,\quad a_{3}=- \frac{g^{2} \left(6 g^{2} k + 5 g^{2} - 2\right)}{2},\\
&a_{4}=- g^{3} \left(2 g^{2} k + 1\right),\quad a_{5}=\frac{g^{6} \left(4 k + 5\right)}{4},\quad a_{6}=g,\quad a_{7}=- \frac{6 g^{2} k + g^{2} + 2}{2},\\
&a_{8}=- g \left(2 g^{2} k + 2 g^{2} - 1\right),\quad a_{9}=\frac{g^{2} \left(4 g^{2} k + g^{2} + 4\right)}{4},\quad a_{10}=g^{3} \left(g - 1\right) \left(g + 1\right).
\end{align*}
\phantomsection\paragraph*{Condition $\CCfam{29}$ (D3-III-08)}\label{center:29}
Free parameters: $g, k$.
Domain: $g\ne0$.
\begin{align*}
&a_{1}=35 g^{2} - 1,\quad a_{2}=- 2 g,\quad a_{3}=- 4 g^{2} \left(336 g^{2} k + 1330 g^{2} - 37\right),\\
&a_{4}=56 g^{3} \left(224 g^{2} k + 360 g^{2} - 5\right),\quad a_{5}=- 112 g^{6} \left(176 k + 175\right),\\
&a_{6}=g,\quad a_{7}=- 48 g^{2} k - 38 g^{2} - 1,\quad a_{8}=- 2 g \left(64 g^{2} k + 4 g^{2} + 1\right),\\
&a_{9}=4 g^{2} \left(1168 g^{2} k + 109 g^{2} + 37\right),\quad a_{10}=- 8 g^{3} \left(1056 g^{2} k + 70 g^{2} + 35\right).
\end{align*}
\phantomsection\paragraph*{Condition $\CCfam{30}$ (D4-I-01)}\label{center:30}
Free parameters: $g, k$.
Domain: $g\ne0$.
\begin{align*}
&a_{1}=- \frac{21 g^{2} + 25}{25},\quad a_{2}=- \frac{18 g}{5},\quad a_{3}=\frac{108 g^{2} \left(224 g^{2} k + 42 g^{2} - 25\right)}{625},\\
&a_{4}=\frac{72 g^{3} \left(7616 g^{2} k + 504 g^{2} - 75\right)}{3125},\quad a_{5}=\frac{9072 g^{6} \left(352 k + 9\right)}{15625},\\
&a_{6}=g,\quad a_{7}=\frac{96 g^{2} k - 54 g^{2} - 25}{25},\quad a_{8}=\frac{2 g \left(896 g^{2} k - 756 g^{2} - 225\right)}{125},\\
&a_{9}=\frac{4 g^{2} \left(992 g^{2} k - 2403 g^{2} - 675\right)}{625},\quad a_{10}=- \frac{72 g^{3} \left(704 g^{2} k + 270 g^{2} + 75\right)}{3125}.
\end{align*}
\phantomsection\paragraph*{Condition $\CCfam{31}$ (D4-I-02)}\label{center:31}
Free parameters: $g, k$.
Domain: $g\ne0$.
\begin{align*}
&a_{1}=- \frac{8 g^{2} + 9}{9},\quad a_{2}=- \frac{8 g}{3},\quad a_{3}=\frac{64 g^{2} \left(12 g^{2} k + 28 g^{2} + 9\right)}{81},\quad a_{4}=\frac{512 g^{3} \left(20 g^{2} k + 20 g^{2} + 9\right)}{243},\\
&a_{5}=\frac{2048 g^{6} \left(17 k + 1\right)}{729},\quad a_{6}=g,\quad a_{7}=\frac{24 g^{2} k - 8 g^{2} - 9}{9},\quad a_{8}=\frac{8 g \left(28 g^{2} k - 20 g^{2} - 9\right)}{27},\\
&a_{9}=- \frac{64 g^{2} \left(3 g^{2} k - 13 g^{2} - 9\right)}{81},\quad a_{10}=- \frac{256 g^{3} \left(17 g^{2} k - 31 g^{2} - 18\right)}{243}.
\end{align*}
\phantomsection\paragraph*{Condition $\CCfam{32}$ (D4-I-05)}\label{center:32}
Free parameters: $g, k$.
Domain: $g\ne0$.
\begin{align*}
&a_{1}=- \frac{13 g^{2} + 49}{49},\quad a_{2}=- \frac{6 g}{7},\quad a_{3}=\frac{12 g^{2} \left(1248 g^{2} k + 26 g^{2} - 49\right)}{2401},\\
&a_{4}=- \frac{8 g^{3} \left(12480 g^{2} k - 104 g^{2} + 49\right)}{16807},\quad a_{5}=\frac{624 g^{6} \left(352 k + 1\right)}{117649},\\
&a_{6}=g,\quad a_{7}=\frac{96 g^{2} k - 6 g^{2} - 49}{49},\quad a_{8}=\frac{2 g \left(640 g^{2} k - 92 g^{2} - 147\right)}{343},\\
&a_{9}=- \frac{4 g^{2} \left(2848 g^{2} k + 115 g^{2} + 147\right)}{2401},\quad a_{10}=\frac{8 g^{3} \left(3520 g^{2} k - 42 g^{2} - 49\right)}{16807}.
\end{align*}
\phantomsection\paragraph*{Condition $\CCfam{33}$ (D4-I-06)}\label{center:33}
Free parameters: $g, k$.
Domain: $g\ne0$.
\begin{align*}
&a_{1}=- \frac{16 g^{2} + 25}{25},\quad a_{2}=0,\quad a_{3}=0,\quad a_{4}=0,\quad a_{5}=0,\quad a_{6}=g,\\
&a_{7}=\frac{24 g^{2} k + 24 g^{2} - 25}{25},\quad a_{8}=\frac{32 g^{3} \left(k + 1\right)}{25},\quad a_{9}=\frac{64 g^{4} \left(k + 1\right)}{125},\quad a_{10}=\frac{256 g^{5} \left(k + 1\right)}{3125}.
\end{align*}
\phantomsection\paragraph*{Condition $\CCfam{34}$ (D4-I-07)}\label{center:34}
Free parameters: $g$.
Domain: $g\ne0$.
\begin{align*}
&a_{1}=- \frac{77 g^{2} + 81}{81},\quad a_{2}=- \frac{10 g}{9},\quad a_{3}=\frac{4 g^{2} \left(32494 g^{2} + 22437\right)}{6561},\quad a_{4}=\frac{8 g^{3} \left(37576 g^{2} + 30483\right)}{19683},\\
&a_{5}=- \frac{80 g^{4} \left(265573 g^{2} + 183384\right)}{531441},\quad a_{6}=g,\quad a_{7}=\frac{2498 g^{2} + 1539}{405},\\
&a_{8}=\frac{2 g \left(6916 g^{2} + 5373\right)}{1215},\quad a_{9}=- \frac{68 g^{2} \left(3311 g^{2} + 1863\right)}{32805},\quad a_{10}=- \frac{8 g^{3} \left(41426 g^{2} + 29133\right)}{19683}.
\end{align*}
\phantomsection\paragraph*{Condition $\CCfam{35}$ (D4-I-09)}\label{center:35}
Free parameters: $g$.
Domain: $g\ne0$.
\begin{align*}
&a_{1}=- \frac{85 g^{2} + 121}{121},\quad a_{2}=\frac{2 g}{11},\quad a_{3}=\frac{4 g^{2} \left(170 g^{2} + 319\right)}{1331},\quad a_{4}=\frac{8 g^{3} \left(15640 g^{2} + 28193\right)}{483153},\\
&a_{5}=\frac{16 g^{4} \left(85 g^{2} + 152\right)}{43923},\quad a_{6}=g,\quad a_{7}=- \frac{3 \left(178 g^{2} + 363\right)}{121},\\
&a_{8}=- \frac{2 g \left(940 g^{2} + 1727\right)}{363},\quad a_{9}=- \frac{4 g^{2} \left(1915 g^{2} + 3443\right)}{3993},\quad a_{10}=- \frac{8 g^{3} \left(4130 g^{2} + 7381\right)}{161051}.
\end{align*}
\phantomsection\paragraph*{Condition $\CCfam{36}$ (D4-II-01)}\label{center:36}
Free parameters: $d, g$.
Domain: $3 d - 11\ne0,\ g\ne0$.
\begin{align*}
&a_{1}=- \frac{21 g^{2} + 25}{25},\quad a_{2}=- \frac{6 g \left(13 d - 5\right)}{5 \left(3 d - 11\right)},\quad a_{3}=\frac{\mathcal N_{\mathrm{a},363}}{1875 d - 6875},\\
&\mathcal N_{\mathrm{a},363}=- 60291 g^{4}- 56400 g^{2}- 44793 d^{2} g^{4}
         - 18900 d^{2} g^{2}+5103 d^{3} g^{4}+77700 d g^{2}+112077 d g^{4},\\
&a_{4}=\frac{\mathcal N_{\mathrm{a},364}}{1875 d - 6875},\\
&\mathcal N_{\mathrm{a},364}=- 12474 g^{5}- 10800 g^{3}- 182574 d^{2} g^{5}- 113400 d^{2} g^{3}
         +30618 d^{3} g^{5}+244200 d g^{3}+261198 d g^{5},\\
&a_{5}=\frac{\mathcal N_{\mathrm{a},365}}{46875 d - 171875},\quad \mathcal N_{\mathrm{a},365}=432 d g^{4} \left(2646 d^{2} g^{2} - 10269 d g^{2} - 9800 d + 2079 g^{2} + 1800\right),\\
&a_{6}=g,\quad a_{7}=- \frac{45 d g^{2} - 207 g^{2} - 175}{25},\quad a_{8}=\frac{\mathcal N_{\mathrm{a},368}}{225 d - 825},\\
&\mathcal N_{\mathrm{a},368}=- 2 g \left(972 d^{2} g^{2} - 5076 d g^{2} - 3735 d + 5544 g^{2} + 5375\right),\quad a_{9}=\frac{\mathcal N_{\mathrm{a},369}}{5625 d - 20625},\\
&\mathcal N_{\mathrm{a},369}=95600 g^{2}+99099 g^{4}- 29403 d^{2} g^{4}- 2187 d^{3} g^{4}
         +8100 d^{2} g^{2}+110187 d g^{4}+146700 d g^{2},\\
&a_{10}=\frac{\mathcal N_{\mathrm{a},3610}}{9375 d - 34375},\\
&\mathcal N_{\mathrm{a},3610}=32400 g^{3}+37422 g^{5}- 251370 d g^{5}- 238200 d g^{3}
         - 10206 d^{3} g^{5}+37800 d^{2} g^{3}+103194 d^{2} g^{5}.
\end{align*}
\phantomsection\paragraph*{Condition $\CCfam{37}$ (D4-III-01)}\label{center:37}
Free parameters: $g, k$.
Domain: $g\ne0$.
\begin{align*}
&a_{1}=\frac{15 g^{2} - 49}{49},\quad a_{2}=\frac{6 g}{7},\quad a_{3}=- \frac{12 g^{2} \left(96 g^{2} k + 30 g^{2} + 49\right)}{2401},\\
&a_{4}=- \frac{8 g^{3} \left(192 g^{2} k - 120 g^{2} - 49\right)}{16807},\quad a_{5}=- \frac{144 g^{6} \left(32 k + 5\right)}{117649},\quad a_{6}=g,\\
&a_{7}=- \frac{96 g^{2} k + 6 g^{2} + 49}{49},\quad a_{8}=- \frac{2 g \left(256 g^{2} k + 76 g^{2} - 147\right)}{343},\\
&a_{9}=- \frac{4 g^{2} \left(224 g^{2} k - 109 g^{2} + 147\right)}{2401},\quad a_{10}=- \frac{8 g^{3} \left(192 g^{2} k + 42 g^{2} - 49\right)}{16807}.
\end{align*}
\phantomsection\paragraph*{Condition $\CCfam{38}$ (D4-III-02)}\label{center:38}
Free parameters: $g, k$.
Domain: $g\ne0$.
\begin{align*}
&a_{1}=\frac{55 g^{2} - 9}{9},\quad a_{2}=- \frac{2 g}{3},\quad a_{3}=- \frac{4 g^{2} \left(1056 g^{2} k + 4730 g^{2} - 765\right)}{81},\\
&a_{4}=\frac{88 g^{3} \left(704 g^{2} k + 1800 g^{2} - 225\right)}{243},\quad a_{5}=- \frac{176 g^{6} \left(1376 k + 1375\right)}{729},\\
&a_{6}=g,\quad a_{7}=- \frac{96 g^{2} k + 86 g^{2} + 9}{9},\quad a_{8}=- \frac{2 g \left(256 g^{2} k + 156 g^{2} + 9\right)}{27},\\
&a_{9}=\frac{4 g^{2} \left(5664 g^{2} k + 1885 g^{2} + 765\right)}{81},\quad a_{10}=- \frac{40 g^{3} \left(2752 g^{2} k + 330 g^{2} + 495\right)}{243}.
\end{align*}
\phantomsection\paragraph*{Condition $\CCfam{39}$ (D5-01)}\label{center:39}
Free parameters: $g$.
Domain: $g\ne0$.
\begin{align*}
&a_{1}=- \frac{45 g^{2} + 49}{49},\quad a_{2}=- \frac{22 g}{7},\quad a_{3}=\frac{12 g^{2} \left(13050 g^{2} + 5047\right)}{2401},\quad a_{4}=\frac{8 g^{3} \left(527400 g^{2} + 244951\right)}{16807},\\
&a_{5}=\frac{16 g^{4} \left(273825 g^{2} + 131768\right)}{16807},\quad a_{6}=g,\quad a_{7}=\frac{3230 g^{2} + 1421}{539},\\
&a_{8}=\frac{2 g \left(4500 g^{2} + 2513\right)}{539},\quad a_{9}=\frac{4 g^{2} \left(17775 g^{2} + 21511\right)}{26411},\quad a_{10}=- \frac{24 g^{3} \left(178650 g^{2} + 80801\right)}{184877}.
\end{align*}
\phantomsection\paragraph*{Condition $\CCfam{40}$ (D5-02)}\label{center:40}
Free parameters: $g$.
Domain: $g\ne0$.
\begin{align*}
&a_{1}=- \frac{65 g^{2} + 81}{81},\quad a_{2}=- \frac{14 g}{9},\quad a_{3}=\frac{4 g^{2} \left(38350 g^{2} + 26811\right)}{6561},\quad a_{4}=\frac{8 g^{3} \left(75400 g^{2} + 56349\right)}{19683},\\
&a_{5}=- \frac{16 g^{4} \left(1284725 g^{2} + 1058346\right)}{531441},\quad a_{6}=g,\quad a_{7}=\frac{3070 g^{2} + 1863}{567},\\
&a_{8}=\frac{2 g \left(9500 g^{2} + 7317\right)}{1701},\quad a_{9}=- \frac{4 g^{2} \left(28925 g^{2} + 21303\right)}{45927},\quad a_{10}=- \frac{8 g^{3} \left(215150 g^{2} + 171909\right)}{137781}.
\end{align*}
\phantomsection\paragraph*{Condition $\CCfam{41}$ (D5-03)}\label{center:41}
Free parameters: $g$.
Domain: $g\ne0$.
\begin{align*}
&a_{1}=- \frac{85 g^{2} + 121}{121},\quad a_{2}=- \frac{6 g}{11},\quad a_{3}=\frac{4 g^{2} \left(26350 g^{2} + 29161\right)}{14641},\\
&a_{4}=- \frac{8 g^{3} \left(98600 g^{2} + 146531\right)}{483153},\quad a_{5}=\frac{16 g^{4} \left(301325 g^{2} + 330572\right)}{5314683},\\
&a_{6}=g,\quad a_{7}=\frac{2270 g^{2} + 2057}{363},\quad a_{8}=\frac{2 g \left(500 g^{2} + 533\right)}{99},\\
&a_{9}=- \frac{4 g^{2} \left(84475 g^{2} + 119911\right)}{131769},\quad a_{10}=\frac{8 g^{3} \left(203150 g^{2} + 232199\right)}{1449459}.
\end{align*}
\phantomsection\paragraph*{Condition $\CCfam{42}$ (D5-04)}\label{center:42}
Free parameters: $g$.
Domain: $g\ne0$.
\begin{align*}
&a_{1}=- \frac{105 g^{2} + 169}{169},\quad a_{2}=\frac{2 g}{13},\quad a_{3}=\frac{12 g^{2} \left(1050 g^{2} + 2197\right)}{28561},\quad a_{4}=\frac{8 g^{3} \left(12600 g^{2} + 25519\right)}{371293},\\
&a_{5}=\frac{16 g^{4} \left(16275 g^{2} + 32786\right)}{4826809},\quad a_{6}=g,\quad a_{7}=- \frac{830 g^{2} + 1859}{169},\quad a_{8}=- \frac{2 g \left(500 g^{2} + 1027\right)}{169},\\
&a_{9}=- \frac{4 g^{2} \left(18475 g^{2} + 37349\right)}{28561},\quad a_{10}=- \frac{8 g^{3} \left(19550 g^{2} + 39377\right)}{371293}.
\end{align*}


\endgroup

\section{Explicit sufficiency certificates}\label{app:suffdata}
Each certificate is paired with the same-numbered condition in \ref{app:atlas}. The row data specify the rational functions used in Section~\ref{sec:sufficiency}; they do not impose further independent conditions on the original coefficients. Identity \eqref{eq:app-low-test} is the verification target for conditions $5$--$13$, and identities \eqref{eq:app-high-defs}--\eqref{eq:app-high-test} are the verification targets for conditions $14$--$42$.
\paragraph*{Degree-two reversibility data ($\CCfam5$--$\CCfam{13}$)}
\begingroup\small
\paragraph*{Certificate for \hyperref[center:5]{$\CCfam{5}$} (D2-swap\_1)}\label{suff:5}
\begin{center}\begin{adjustbox}{max width=\linewidth}
$\displaystyle\begin{aligned}
&b_{5}=4 g,\qquad F_{5}(s)=- \frac{\left(4 g s + 1\right)^{4} \left(6 g s - 1\right)}{10 g s + 1},\qquad \rho_{5}(s)=1.
\end{aligned}$
\end{adjustbox}\end{center}
\paragraph*{Certificate for \hyperref[center:6]{$\CCfam{6}$} (D2-swap\_2)}\label{suff:6}
\begin{center}\begin{adjustbox}{max width=\linewidth}
$\displaystyle\begin{aligned}
&b_{6}=4 g,\qquad F_{6}(s)=\frac{\left(2 g s - 1\right)^{2} \left(4 g s + 1\right)^{4}}{\left(6 g s + 1\right)^{2}},\qquad \rho_{6}(s)=1.
\end{aligned}$
\end{adjustbox}\end{center}
\paragraph*{Certificate for \hyperref[center:7]{$\CCfam{7}$} (D2-swap\_3)}\label{suff:7}
\begin{center}\begin{adjustbox}{max width=\linewidth}
$\displaystyle\begin{aligned}
&b_{7}=4 g,\qquad F_{7}(s)=- \frac{\left(2 g s - 3\right)^{3} \left(4 g s + 1\right)^{4}}{\left(14 g s + 3\right)^{3}},\qquad \rho_{7}(s)=1.
\end{aligned}$
\end{adjustbox}\end{center}
\paragraph*{Certificate for \hyperref[center:8]{$\CCfam{8}$} (D2-fixed\_1)}\label{suff:8}
\begin{center}\begin{adjustbox}{max width=\linewidth}
$\displaystyle\begin{aligned}
&b_{8}=- 4 g,\qquad F_{8}(s)=- \frac{\left(4 g s - 1\right)^{4} \left(6 g s + 1\right)}{10 g s - 1},\qquad \rho_{8}(s)=\frac{- 6 g s - 1}{10 g s - 1}.
\end{aligned}$
\end{adjustbox}\end{center}
\paragraph*{Certificate for \hyperref[center:9]{$\CCfam{9}$} (D2-fixed\_3)}\label{suff:9}
\begin{center}\begin{adjustbox}{max width=\linewidth}
$\displaystyle\begin{aligned}
&b_{9}=\frac{4 g}{3},\qquad F_{9}(s)=\frac{\left(- 2 g s + 1\right) \left(4 g s + 3\right)^{4}}{270 g s + 81},\qquad \rho_{9}(s)=\frac{- 10 g s - 3}{6 g s - 3}.
\end{aligned}$
\end{adjustbox}\end{center}
\paragraph*{Certificate for \hyperref[center:10]{$\CCfam{10}$} (D2-fixed\_4)}\label{suff:10}
\begin{center}\begin{adjustbox}{max width=\linewidth}
$\displaystyle\begin{aligned}
&b_{10}=2 g,\qquad F_{10}(s)=\frac{\left(g s - 1\right)^{2} \left(2 g s + 1\right)^{4}}{\left(3 g s + 1\right)^{2}},\qquad \rho_{10}(s)=\frac{- 3 g s - 1}{g s - 1}.
\end{aligned}$
\end{adjustbox}\end{center}
\paragraph*{Certificate for \hyperref[center:11]{$\CCfam{11}$} (D2-fixed\_5)}\label{suff:11}
\begin{center}\begin{adjustbox}{max width=\linewidth}
$\displaystyle\begin{aligned}
&b_{11}=\frac{12 g}{5},\qquad F_{11}(s)=- \frac{\left(2 g s - 5\right)^{3} \left(12 g s + 5\right)^{4}}{625 \left(14 g s + 5\right)^{3}},\qquad \rho_{11}(s)=\frac{- 14 g s - 5}{2 g s - 5}.
\end{aligned}$
\end{adjustbox}\end{center}
\paragraph*{Certificate for \hyperref[center:12]{$\CCfam{12}$} (D2-fixed\_6)}\label{suff:12}
\begin{center}\begin{adjustbox}{max width=\linewidth}
$\displaystyle\begin{aligned}
&b_{12}=\frac{8 g}{3},\qquad F_{12}(s)=1,\qquad \rho_{12}(s)=\frac{8 g s}{3} + 1.
\end{aligned}$
\end{adjustbox}\end{center}
\paragraph*{Certificate for \hyperref[center:13]{$\CCfam{13}$} (D2-fixed\_7)}\label{suff:13}
\begin{center}\begin{adjustbox}{max width=\linewidth}
$\displaystyle\begin{aligned}
&b_{13}=\frac{20 g}{7},\qquad F_{13}(s)=\frac{\left(2 g s + 7\right)^{5} \left(20 g s + 7\right)^{4}}{2401 \left(18 g s + 7\right)^{5}},\qquad \rho_{13}(s)=\frac{18 g s + 7}{2 g s + 7}.
\end{aligned}$
\end{adjustbox}\end{center}
\endgroup
\paragraph*{Higher-degree rational pullback data ($\CCfam{14}$--$\CCfam{42}$)}
\begingroup\small
\setlength{\jot}{2pt}
\allowdisplaybreaks[4]
\setlength{\abovedisplayskip}{4pt}
\setlength{\belowdisplayskip}{4pt}
\par\addvspace{6pt}\phantomsection\noindent\textit{Certificate for \hyperref[center:14]{$\CCfam{14}$} (D3-I-02).}\label{suff:14}\nopagebreak[4]

\begin{align*}
&c_{14}=- \frac{3 g^{2}}{4},\qquad \alpha_{14}=- g,\\
&\Omega_{14}=- 24 g^{2},\qquad \lambda_{14}=1,\\
&\kappa_{14}=k,\qquad \nu_{14}=u,\\
&\gamma_{14}=- \frac{1}{3 g^{2}},\\
&E_{14}(z)=\left(z - 1\right) \left(3 z - 1\right) \left(3 z + 1\right),\\
&N_{14}(z)=\left(3 z - 1\right)^{2} \left(3 z + 1\right),\\
&V_{14}(z)=- \left(z - 1\right)^{3}.
\end{align*}
\par\addvspace{6pt}\phantomsection\noindent\textit{Certificate for \hyperref[center:15]{$\CCfam{15}$} (D3-I-03).}\label{suff:15}\nopagebreak[4]

\begin{align*}
&c_{15}=- \frac{21 g^{2}}{25},\qquad \alpha_{15}=- \frac{6 g}{5},\\
&\Omega_{15}=- \frac{384 g^{2}}{25},\qquad \lambda_{15}=1,\\
&\kappa_{15}=k,\qquad \nu_{15}=- \frac{36 g^{2} k + 36 g^{2} + 25}{48 g^{2}},\\
&\gamma_{15}=- \frac{25}{48 g^{2}},\\
&E_{15}(z)=\frac{\left(z - 1\right)^{2} \left(5 z + 3\right)}{3},\\
&N_{15}(z)=\frac{\left(5 z + 3\right) \left(7 z - 3\right)^{2}}{27},\\
&V_{15}(z)=- \left(z - 1\right)^{3}.
\end{align*}
\par\addvspace{6pt}\phantomsection\noindent\textit{Certificate for \hyperref[center:16]{$\CCfam{16}$} (D3-I-04).}\label{suff:16}\nopagebreak[4]

\begin{align*}
&c_{16}=\frac{7 g^{2}}{9},\qquad \alpha_{16}=\frac{2 g}{3},\\
&\Omega_{16}=- \frac{128 g^{2}}{3},\qquad \lambda_{16}=\frac{1}{2},\\
&\kappa_{16}=k,\qquad \nu_{16}=- \frac{7 k}{9},\\
&\gamma_{16}=- \frac{1}{4 g^{2}},\\
&E_{16}(z)=- \left(5 z - 1\right) \left(7 z + 1\right)^{2},\\
&N_{16}(z)=\left(5 z - 1\right)^{2} \left(7 z + 1\right),\\
&V_{16}(z)=- \left(z - 1\right)^{3}.
\end{align*}
\par\addvspace{6pt}\phantomsection\noindent\textit{Certificate for \hyperref[center:17]{$\CCfam{17}$} (D3-I-05).}\label{suff:17}\nopagebreak[4]

\begin{align*}
&c_{17}=- \frac{9 g^{2}}{25},\qquad \alpha_{17}=- \frac{2 g}{5},\\
&\Omega_{17}=- \frac{384 g^{2}}{25},\qquad \lambda_{17}=\frac{1}{2},\\
&\kappa_{17}=k,\qquad \nu_{17}=u,\\
&\gamma_{17}=- \frac{25}{36 g^{2}},\\
&E_{17}(z)=\left(z - 1\right) \left(3 z + 1\right) \left(9 z - 1\right),\\
&N_{17}(z)=- \left(3 z + 1\right)^{2} \left(9 z - 1\right),\\
&V_{17}(z)=- \left(z - 1\right)^{3}.
\end{align*}
\par\addvspace{6pt}\phantomsection\noindent\textit{Certificate for \hyperref[center:18]{$\CCfam{18}$} (D3-I-06).}\label{suff:18}\nopagebreak[4]

\begin{align*}
&c_{18}=- \frac{33 g^{2}}{49},\qquad \alpha_{18}=- \frac{6 g}{7},\\
&\Omega_{18}=- \frac{384 g^{2}}{49},\qquad \lambda_{18}=\frac{1}{2},\\
&\kappa_{18}=k,\qquad \nu_{18}=- \frac{36 g^{2} k + 36 g^{2} + 49}{36 g^{2}},\\
&\gamma_{18}=- \frac{49}{36 g^{2}},\\
&E_{18}(z)=\frac{\left(z - 1\right)^{2} \left(z + 3\right)}{3},\\
&N_{18}(z)=- \frac{\left(z + 3\right)^{2} \left(11 z - 3\right)}{27},\\
&V_{18}(z)=- \left(z - 1\right)^{3}.
\end{align*}
\par\addvspace{6pt}\phantomsection\noindent\textit{Certificate for \hyperref[center:19]{$\CCfam{19}$} (D3-I-07).}\label{suff:19}\nopagebreak[4]

\begin{align*}
&c_{19}=- \frac{8 g^{2}}{9},\qquad \alpha_{19}=- \frac{4 g}{3},\\
&\Omega_{19}=- \frac{32 g^{2}}{3},\qquad \lambda_{19}=2,\\
&\kappa_{19}=k,\qquad \nu_{19}=- \frac{8 k}{9},\\
&\gamma_{19}=- \frac{1}{2 g^{2}},\\
&E_{19}(z)=- \left(z + 1\right)^{2} \left(2 z - 1\right),\\
&N_{19}(z)=\left(z + 1\right) \left(2 z - 1\right)^{2},\\
&V_{19}(z)=- \left(z - 1\right)^{3}.
\end{align*}
\par\addvspace{6pt}\phantomsection\noindent\textit{Certificate for \hyperref[center:20]{$\CCfam{20}$} (D3-I-08).}\label{suff:20}\nopagebreak[4]

\begin{align*}
&c_{20}=- \frac{45 g^{2}}{49},\qquad \alpha_{20}=- \frac{10 g}{7},\\
&\Omega_{20}=- \frac{384 g^{2}}{49},\qquad \lambda_{20}=2,\\
&\kappa_{20}=k,\qquad \nu_{20}=- \frac{60 g^{2} k + 60 g^{2} + 49}{72 g^{2}},\\
&\gamma_{20}=- \frac{49}{72 g^{2}},\\
&E_{20}(z)=- \frac{\left(z - 1\right) \left(3 z + 5\right)^{2}}{25},\\
&N_{20}(z)=\frac{\left(3 z + 5\right) \left(9 z - 5\right)^{2}}{125},\\
&V_{20}(z)=- \left(z - 1\right)^{3}.
\end{align*}
\par\addvspace{6pt}\phantomsection\noindent\textit{Certificate for \hyperref[center:21]{$\CCfam{21}$} (D3-I-11).}\label{suff:21}\nopagebreak[4]

\begin{align*}
&c_{21}=- \frac{77 g^{2}}{81},\qquad \alpha_{21}=- \frac{14 g}{9},\\
&\Omega_{21}=- \frac{128 g^{2}}{27},\qquad \lambda_{21}=3,\\
&\kappa_{21}=\frac{3 \left(28 g^{2} + 27\right)}{4 g^{2}},\\
&\nu_{21}=- \frac{11 \left(28 g^{2} + 27\right)}{16 g^{2}},\qquad \gamma_{21}=- \frac{27}{32 g^{2}},\\
&E_{21}(z)=\frac{\left(z + 7\right)^{3}}{343},\\
&N_{21}(z)=\frac{\left(z + 7\right) \left(11 z - 7\right)^{2}}{343},\\
&V_{21}(z)=- \left(z - 1\right)^{3}.
\end{align*}
\par\addvspace{6pt}\phantomsection\noindent\textit{Certificate for \hyperref[center:22]{$\CCfam{22}$} (D3-I-12).}\label{suff:22}\nopagebreak[4]

\begin{align*}
&c_{22}=- \frac{65 g^{2}}{81},\qquad \alpha_{22}=- \frac{10 g}{9},\\
&\Omega_{22}=- \frac{128 g^{2}}{27},\qquad \lambda_{22}=\frac{3}{2},\\
&\kappa_{22}=- \frac{3 \left(20 g^{2} + 27\right)}{8 g^{2}},\\
&\nu_{22}=\frac{13 \left(20 g^{2} + 27\right)}{40 g^{2}},\qquad \gamma_{22}=- \frac{27}{20 g^{2}},\\
&E_{22}(z)=- \frac{\left(z - 5\right)^{3}}{125},\\
&N_{22}(z)=- \frac{\left(z - 5\right)^{2} \left(13 z - 5\right)}{125},\\
&V_{22}(z)=- \left(z - 1\right)^{3}.
\end{align*}
\par\addvspace{6pt}\phantomsection\noindent\textit{Certificate for \hyperref[center:23]{$\CCfam{23}$} (D3-II-01).}\label{suff:23}\nopagebreak[4]

\begin{align*}
&c_{23}=- \frac{5 g^{2}}{9},\qquad \alpha_{23}=- \frac{2 g}{3},\\
&\Omega_{23}=\frac{8 g^{2} \left(d - 13\right)}{9},\qquad \lambda_{23}=1,\\
&\kappa_{23}=k,\qquad \nu_{23}=- \frac{5 k}{6},\\
&\gamma_{23}=\frac{9}{g^{2} \left(d - 13\right)},\\
&E_{23}(z)=- \left(5 z - 1\right) \left(d z^{2} + 2 z + 1\right),\\
&N_{23}(z)=- \left(5 z - 1\right) \left(d z^{2} + 2 z + 1\right),\\
&V_{23}(z)=- \left(z - 1\right)^{3}.
\end{align*}
\par\addvspace{6pt}\phantomsection\noindent\textit{Certificate for \hyperref[center:24]{$\CCfam{24}$} (D3-II-02).}\label{suff:24}\nopagebreak[4]

\begin{align*}
&c_{24}=- \frac{3 g^{2}}{4},\qquad \alpha_{24}=- g,\\
&\Omega_{24}=2 g^{2} \left(d - 3\right),\qquad \lambda_{24}=1,\\
&\kappa_{24}=k,\\
&\nu_{24}=- \frac{d g^{2} k + d g^{2} - 3 g^{2} k - 3 g^{2} - 4}{g^{2} \left(d - 3\right)},\\
&\gamma_{24}=\frac{4}{g^{2} \left(d - 3\right)},\\
&E_{24}(z)=- \left(z - 1\right) \left(d z^{2} + 1\right),\\
&N_{24}(z)=- \left(3 z - 1\right) \left(d z^{2} + 1\right),\\
&V_{24}(z)=- \left(z - 1\right)^{3}.
\end{align*}
\par\addvspace{6pt}\phantomsection\noindent\textit{Certificate for \hyperref[center:25]{$\CCfam{25}$} (D3-III-01).}\label{suff:25}\nopagebreak[4]

\begin{align*}
&c_{25}=\frac{11 g^{2}}{25},\qquad \alpha_{25}=- \frac{22 g}{5},\\
&\Omega_{25}=\frac{384 g^{2}}{25},\qquad \lambda_{25}=\frac{1}{3},\\
&\kappa_{25}=k,\qquad \nu_{25}=- \frac{k}{4},\\
&\gamma_{25}=\frac{25}{32 g^{2}},\\
&E_{25}(z)=- \frac{\left(z + 11\right)^{2} \left(3 z - 11\right)}{1331},\\
&N_{25}(z)=- \frac{\left(3 z - 11\right)^{3}}{1331},\\
&V_{25}(z)=- \frac{\left(z - 1\right) \left(z + 11\right)^{2}}{121}.
\end{align*}
\par\addvspace{6pt}\phantomsection\noindent\textit{Certificate for \hyperref[center:26]{$\CCfam{26}$} (D3-III-02).}\label{suff:26}\nopagebreak[4]

\begin{align*}
&c_{26}=8 g^{2},\qquad \alpha_{26}=- 8 g,\\
&\Omega_{26}=96 g^{2},\qquad \lambda_{26}=\frac{1}{3},\\
&\kappa_{26}=k,\qquad \nu_{26}=u,\\
&\gamma_{26}=\frac{1}{8 g^{2}},\\
&E_{26}(z)=- \frac{\left(z - 1\right) \left(z + 2\right)}{2},\\
&N_{26}(z)=1,\\
&V_{26}(z)=- \frac{\left(z - 1\right) \left(z + 2\right)^{2}}{4}.
\end{align*}
\par\addvspace{6pt}\phantomsection\noindent\textit{Certificate for \hyperref[center:27]{$\CCfam{27}$} (D3-III-03).}\label{suff:27}\nopagebreak[4]

\begin{align*}
&c_{27}=35 g^{2},\qquad \alpha_{27}=10 g,\\
&\Omega_{27}=384 g^{2},\qquad \lambda_{27}=\frac{1}{3},\\
&\kappa_{27}=k,\qquad \nu_{27}=- \frac{5 k}{8},\\
&\gamma_{27}=\frac{1}{32 g^{2}},\\
&E_{27}(z)=\frac{\left(z - 1\right)^{2} \left(3 z + 5\right)}{5},\\
&N_{27}(z)=\frac{\left(3 z + 5\right)^{3}}{125},\\
&V_{27}(z)=- \frac{\left(z - 1\right) \left(7 z + 5\right)^{2}}{25}.
\end{align*}
\par\addvspace{6pt}\phantomsection\noindent\textit{Certificate for \hyperref[center:28]{$\CCfam{28}$} (D3-III-07).}\label{suff:28}\nopagebreak[4]

\begin{align*}
&c_{28}=\frac{5 g^{2}}{4},\qquad \alpha_{28}=- 5 g,\\
&\Omega_{28}=24 g^{2},\qquad \lambda_{28}=\frac{2}{3},\\
&\kappa_{28}=k,\qquad \nu_{28}=- \frac{2 k}{5},\\
&\gamma_{28}=\frac{2}{5 g^{2}},\\
&E_{28}(z)=\frac{\left(z - 5\right)^{2} \left(z + 5\right)}{125},\\
&N_{28}(z)=- \frac{\left(z - 5\right)^{3}}{125},\\
&V_{28}(z)=- \frac{\left(z - 1\right) \left(z + 5\right)^{2}}{25}.
\end{align*}
\par\addvspace{6pt}\phantomsection\noindent\textit{Certificate for \hyperref[center:29]{$\CCfam{29}$} (D3-III-08).}\label{suff:29}\nopagebreak[4]

\begin{align*}
&c_{29}=35 g^{2},\qquad \alpha_{29}=- 14 g,\\
&\Omega_{29}=384 g^{2},\qquad \lambda_{29}=\frac{2}{3},\\
&\kappa_{29}=k,\qquad \nu_{29}=- \frac{7 k}{10},\\
&\gamma_{29}=\frac{1}{40 g^{2}},\\
&E_{29}(z)=- \frac{\left(z - 1\right) \left(z + 7\right)^{2}}{49},\\
&N_{29}(z)=\frac{\left(z + 7\right)^{3}}{343},\\
&V_{29}(z)=- \frac{\left(z - 1\right) \left(5 z + 7\right)^{2}}{49}.
\end{align*}
\par\addvspace{6pt}\phantomsection\noindent\textit{Certificate for \hyperref[center:30]{$\CCfam{30}$} (D4-I-01).}\label{suff:30}\nopagebreak[4]

\begin{align*}
&c_{30}=- \frac{21 g^{2}}{25},\qquad \alpha_{30}=- \frac{6 g}{5},\\
&\Omega_{30}=- \frac{768 g^{2}}{25},\qquad \lambda_{30}=1,\\
&\kappa_{30}=k,\qquad \nu_{30}=- \frac{7 k}{8},\\
&\gamma_{30}=- \frac{25}{96 g^{2}},\\
&E_{30}(z)=- \frac{\left(3 z + 1\right) \left(7 z - 3\right)}{3},\\
&N_{30}(z)=- \frac{\left(3 z + 1\right) \left(7 z - 3\right)^{3}}{27},\\
&V_{30}(z)=\left(z - 1\right)^{4}.
\end{align*}
\par\addvspace{6pt}\phantomsection\noindent\textit{Certificate for \hyperref[center:31]{$\CCfam{31}$} (D4-I-02).}\label{suff:31}\nopagebreak[4]

\begin{align*}
&c_{31}=- \frac{8 g^{2}}{9},\qquad \alpha_{31}=- \frac{4 g}{3},\\
&\Omega_{31}=- \frac{64 g^{2}}{3},\qquad \lambda_{31}=1,\\
&\kappa_{31}=k,\qquad \nu_{31}=- \frac{16 g^{2} k + 16 g^{2} + 9}{24 g^{2}},\\
&\gamma_{31}=- \frac{3}{8 g^{2}},\\
&E_{31}(z)=- \left(z - 1\right) \left(2 z + 1\right),\\
&N_{31}(z)=- \left(2 z - 1\right)^{3} \left(2 z + 1\right),\\
&V_{31}(z)=\left(z - 1\right)^{4}.
\end{align*}
\par\addvspace{6pt}\phantomsection\noindent\textit{Certificate for \hyperref[center:32]{$\CCfam{32}$} (D4-I-05).}\label{suff:32}\nopagebreak[4]

\begin{align*}
&c_{32}=- \frac{13 g^{2}}{49},\qquad \alpha_{32}=- \frac{2 g}{7},\\
&\Omega_{32}=- \frac{768 g^{2}}{49},\qquad \lambda_{32}=\frac{1}{3},\\
&\kappa_{32}=k,\qquad \nu_{32}=- \frac{13 k}{16},\\
&\gamma_{32}=- \frac{49}{64 g^{2}},\\
&E_{32}(z)=- \left(3 z + 1\right) \left(13 z - 1\right),\\
&N_{32}(z)=- \left(3 z + 1\right)^{3} \left(13 z - 1\right),\\
&V_{32}(z)=\left(z - 1\right)^{4}.
\end{align*}
\par\addvspace{6pt}\phantomsection\noindent\textit{Certificate for \hyperref[center:33]{$\CCfam{33}$} (D4-I-06).}\label{suff:33}\nopagebreak[4]

\begin{align*}
&c_{33}=- \frac{16 g^{2}}{25},\qquad \alpha_{33}=- \frac{4 g}{5},\\
&\Omega_{33}=- \frac{192 g^{2}}{25},\qquad \lambda_{33}=\frac{1}{3},\\
&\kappa_{33}=k,\qquad \nu_{33}=- \frac{16 g^{2} k + 16 g^{2} + 25}{16 g^{2}},\\
&\gamma_{33}=- \frac{25}{16 g^{2}},\\
&E_{33}(z)=1 - z,\\
&N_{33}(z)=1 - 4 z,\\
&V_{33}(z)=\left(z - 1\right)^{4}.
\end{align*}
\par\addvspace{6pt}\phantomsection\noindent\textit{Certificate for \hyperref[center:34]{$\CCfam{34}$} (D4-I-07).}\label{suff:34}\nopagebreak[4]

\begin{align*}
&c_{34}=- \frac{77 g^{2}}{81},\qquad \alpha_{34}=- \frac{14 g}{9},\\
&\Omega_{34}=- \frac{256 g^{2}}{27},\qquad \lambda_{34}=2,\\
&\kappa_{34}=\frac{112 g^{2} + 81}{20 g^{2}},\\
&\nu_{34}=- \frac{11 \left(112 g^{2} + 81\right)}{240 g^{2}},\qquad \gamma_{34}=- \frac{9}{16 g^{2}},\\
&E_{34}(z)=\frac{\left(5 z + 7\right)^{2}}{49},\\
&N_{34}(z)=- \frac{\left(5 z + 7\right) \left(11 z - 7\right)^{3}}{2401},\\
&V_{34}(z)=\left(z - 1\right)^{4}.
\end{align*}
\par\addvspace{6pt}\phantomsection\noindent\textit{Certificate for \hyperref[center:35]{$\CCfam{35}$} (D4-I-09).}\label{suff:35}\nopagebreak[4]

\begin{align*}
&c_{35}=- \frac{85 g^{2}}{121},\qquad \alpha_{35}=- \frac{10 g}{11},\\
&\Omega_{35}=- \frac{768 g^{2}}{121},\qquad \lambda_{35}=\frac{2}{3},\\
&\kappa_{35}=- \frac{80 g^{2} + 121}{12 g^{2}},\\
&\nu_{35}=\frac{17 \left(80 g^{2} + 121\right)}{240 g^{2}},\qquad \gamma_{35}=- \frac{121}{80 g^{2}},\\
&E_{35}(z)=\frac{\left(z - 5\right)^{2}}{25},\\
&N_{35}(z)=\frac{\left(z - 5\right)^{3} \left(17 z - 5\right)}{625},\\
&V_{35}(z)=\left(z - 1\right)^{4}.
\end{align*}
\par\addvspace{6pt}\phantomsection\noindent\textit{Certificate for \hyperref[center:36]{$\CCfam{36}$} (D4-II-01).}\label{suff:36}\nopagebreak[4]

\begin{align*}
&c_{36}=- \frac{21 g^{2}}{25},\qquad \alpha_{36}=- \frac{6 g}{5},\\
&\Omega_{36}=\frac{24 g^{2} \left(3 d - 11\right)}{25},\qquad \lambda_{36}=1,\\
&\kappa_{36}=\frac{2 \left(27 d g^{2} - 99 g^{2} - 100\right)}{3 g^{2} \left(3 d - 11\right)},\\
&\nu_{36}=- \frac{7 \left(27 d g^{2} - 99 g^{2} - 100\right)}{12 g^{2} \left(3 d - 11\right)},\\
&\gamma_{36}=\frac{25}{3 g^{2} \left(3 d - 11\right)},\\
&E_{36}(z)=\frac{3 d z^{2} + 2 z + 3}{3},\\
&N_{36}(z)=\frac{\left(7 z - 3\right)^{2} \left(3 d z^{2} + 2 z + 3\right)}{27},\\
&V_{36}(z)=\left(z - 1\right)^{4}.
\end{align*}
\par\addvspace{6pt}\phantomsection\noindent\textit{Certificate for \hyperref[center:37]{$\CCfam{37}$} (D4-III-01).}\label{suff:37}\nopagebreak[4]

\begin{align*}
&c_{37}=\frac{15 g^{2}}{49},\qquad \alpha_{37}=- \frac{30 g}{7},\\
&\Omega_{37}=\frac{768 g^{2}}{49},\qquad \lambda_{37}=\frac{1}{4},\\
&\kappa_{37}=k,\qquad \nu_{37}=- \frac{k}{5},\\
&\gamma_{37}=\frac{49}{60 g^{2}},\\
&E_{37}(z)=- \frac{\left(z - 5\right) \left(z + 15\right)}{75},\\
&N_{37}(z)=\frac{\left(z - 5\right)^{4}}{625},\\
&V_{37}(z)=- \frac{\left(z - 1\right) \left(z + 15\right)^{3}}{3375}.
\end{align*}
\par\addvspace{6pt}\phantomsection\noindent\textit{Certificate for \hyperref[center:38]{$\CCfam{38}$} (D4-III-02).}\label{suff:38}\nopagebreak[4]

\begin{align*}
&c_{38}=\frac{55 g^{2}}{9},\qquad \alpha_{38}=- \frac{22 g}{3},\\
&\Omega_{38}=\frac{256 g^{2}}{3},\qquad \lambda_{38}=\frac{1}{4},\\
&\kappa_{38}=k,\qquad \nu_{38}=- \frac{11 k}{15},\\
&\gamma_{38}=\frac{3}{20 g^{2}},\\
&E_{38}(z)=- \frac{\left(z - 1\right) \left(z + 11\right)}{11},\\
&N_{38}(z)=\frac{\left(z + 11\right)^{4}}{14641},\\
&V_{38}(z)=- \frac{\left(z - 1\right) \left(5 z + 11\right)^{3}}{1331}.
\end{align*}
\par\addvspace{6pt}\phantomsection\noindent\textit{Certificate for \hyperref[center:39]{$\CCfam{39}$} (D5-01).}\label{suff:39}\nopagebreak[4]

\begin{align*}
&c_{39}=- \frac{45 g^{2}}{49},\qquad \alpha_{39}=- \frac{10 g}{7},\\
&\Omega_{39}=- \frac{1280 g^{2}}{49},\qquad \lambda_{39}=1,\\
&\kappa_{39}=\frac{100 g^{2} + 49}{44 g^{2}},\\
&\nu_{39}=- \frac{9 \left(100 g^{2} + 49\right)}{440 g^{2}},\qquad \gamma_{39}=- \frac{49}{160 g^{2}},\\
&E_{39}(z)=\frac{11 z + 5}{5},\\
&N_{39}(z)=\frac{\left(9 z - 5\right)^{4} \left(11 z + 5\right)}{3125},\\
&V_{39}(z)=- \left(z - 1\right)^{5}.
\end{align*}
\par\addvspace{6pt}\phantomsection\noindent\textit{Certificate for \hyperref[center:40]{$\CCfam{40}$} (D5-02).}\label{suff:40}\nopagebreak[4]

\begin{align*}
&c_{40}=- \frac{65 g^{2}}{81},\qquad \alpha_{40}=- \frac{10 g}{9},\\
&\Omega_{40}=- \frac{640 g^{2}}{27},\qquad \lambda_{40}=\frac{1}{2},\\
&\kappa_{40}=\frac{100 g^{2} + 81}{56 g^{2}},\\
&\nu_{40}=- \frac{13 \left(100 g^{2} + 81\right)}{840 g^{2}},\qquad \gamma_{40}=- \frac{9}{20 g^{2}},\\
&E_{40}(z)=\frac{7 z + 5}{5},\\
&N_{40}(z)=- \frac{\left(7 z + 5\right)^{2} \left(13 z - 5\right)^{3}}{3125},\\
&V_{40}(z)=- \left(z - 1\right)^{5}.
\end{align*}
\par\addvspace{6pt}\phantomsection\noindent\textit{Certificate for \hyperref[center:41]{$\CCfam{41}$} (D5-03).}\label{suff:41}\nopagebreak[4]

\begin{align*}
&c_{41}=- \frac{85 g^{2}}{121},\qquad \alpha_{41}=- \frac{10 g}{11},\\
&\Omega_{41}=- \frac{1920 g^{2}}{121},\qquad \lambda_{41}=\frac{1}{3},\\
&\kappa_{41}=\frac{100 g^{2} + 121}{36 g^{2}},\\
&\nu_{41}=- \frac{17 \left(100 g^{2} + 121\right)}{720 g^{2}},\\
&\gamma_{41}=- \frac{121}{160 g^{2}},\\
&E_{41}(z)=\frac{3 z + 5}{5},\\
&N_{41}(z)=\frac{\left(3 z + 5\right)^{3} \left(17 z - 5\right)^{2}}{3125},\\
&V_{41}(z)=- \left(z - 1\right)^{5}.
\end{align*}
\par\addvspace{6pt}\phantomsection\noindent\textit{Certificate for \hyperref[center:42]{$\CCfam{42}$} (D5-04).}\label{suff:42}\nopagebreak[4]

\begin{align*}
&c_{42}=- \frac{105 g^{2}}{169},\qquad \alpha_{42}=- \frac{10 g}{13},\\
&\Omega_{42}=- \frac{1280 g^{2}}{169},\qquad \lambda_{42}=\frac{1}{4},\\
&\kappa_{42}=- \frac{100 g^{2} + 169}{16 g^{2}},\\
&\nu_{42}=\frac{21 \left(100 g^{2} + 169\right)}{400 g^{2}},\qquad \gamma_{42}=- \frac{169}{100 g^{2}},\\
&E_{42}(z)=- \frac{z - 5}{5},\\
&N_{42}(z)=- \frac{\left(z - 5\right)^{4} \left(21 z - 5\right)}{3125},\\
&V_{42}(z)=- \left(z - 1\right)^{5}.
\end{align*}
\endgroup

\begingroup
\normalsize
\raggedbottom
\renewcommand{\arraystretch}{1.02}
\setlength{\jot}{2pt}
\setlength{\abovedisplayskip}{4pt}
\setlength{\belowdisplayskip}{4pt}
\setlength{\abovedisplayshortskip}{3pt}
\setlength{\belowdisplayshortskip}{3pt}
\setlength{\LTpre}{5pt}
\setlength{\LTpost}{5pt}
\allowdisplaybreaks[2]
\section{Four-degree records and their destinations}\label{app:four}
For convenience, the degree-four source records are indexed here in one place. The integer $k$ in a record is a discrete mass label and should not be confused with any free real $k$ in an original-coefficient parametrization.
{\small
\begin{longtable}{@{}lllll@{}}
\caption{Degree-four source labels and their retained destinations.}\label{tab:degree4}\\
\toprule Record & $k$ or $j$ & $e$ & $(m,n)$ & Target\\\midrule\endfirsthead
\toprule Record & $k$ or $j$ & $e$ & $(m,n)$ & Target\\\midrule\endhead
D4-I-01 & $k=1$ & 1 & $(0,1)$ & $\CCfam{30}$\\
D4-I-02 & $k=1$ & 1 & $(1,0)$ & $\CCfam{31}$\\
D4-I-03 & $k=1$ & 2 & $(0,1)$ & $\CCfam9$\\
D4-I-04 & $k=1$ & 2 & $(1,0)$ & $\CCfam{10}$\\
D4-I-05 & $k=1$ & 3 & $(0,1)$ & $\CCfam{32}$\\
D4-I-06 & $k=1$ & 3 & $(1,0)$ & $\CCfam{33}$\\
D4-I-07 & $k=2$ & 1 & $(0,0)$ & $\CCfam{34}$\\
D4-I-08 & $k=2$ & 2 & $(0,0)$ & $\CCfam{11}$\\
D4-I-09 & $k=2$ & 3 & $(0,0)$ & $\CCfam{35}$\\
D4-II-01 & -- & -- & -- & $\CCfam{36}$\\
D4-III-01 & $j=1$ & -- & $(0,1)$ & $\CCfam{37}$\\
D4-III-02 & $j=1$ & -- & $(1,0)$ & $\CCfam{38}$\\
D4-III-03 & $j=2$ & -- & $(0,1)$ & $\CCfam5$\\
D4-III-04 & $j=2$ & -- & $(1,0)$ & $\CCfam5$\\
D4-III-05 & $j=3$ & -- & $(0,1)$ & $\CCfam{38}$\\
D4-III-06 & $j=3$ & -- & $(1,0)$ & $\CCfam{37}$\\
D4-B-01 & -- & 1 & -- & $\CCfam{31}$\\
D4-B-02 & -- & 2 & -- & $\CCfam{10}$\\
D4-B-03 & -- & 3 & -- & $\CCfam{33}$\\
\bottomrule
\end{longtable}
}
In Table~\ref{tab:degree4}, a dash indicates that the corresponding indexing parameter is not used in the displayed record; it does not represent zero or an arbitrary value. For Type-III records, the two local degrees in the positive fiber are $j$ and $4-j$, so no separate $e$ is listed. Type-II records and boundary records use different indexing conventions; their omitted entries are determined by the corresponding fiber configurations discussed in the classification proof.

\section{Complete source-record destinations and coefficient transports}\label{app:records}
This index is part of the necessity proof. The source tuple records the integer mass data, not the free parameters of the final rational chart. The arrows refer to the original ten coefficients. A transported source may specialize the target parameters; it need not be identical to the whole target family. The ordered stratum argument in Section~\ref{sec:necessity-order} provides the domains of the four elementary sources and the nine degree-two sources.
{\small
\begin{align*}
&\mathrm{F}:\; \text{proportional}\;\longmapsto\CCfam{1}\quad \mathrm{A}:\; \text{reflection}\;\longmapsto\CCfam{2}\quad \mathrm{G}:\; \text{composition}\;\longmapsto\CCfam{3}\\
&\mathrm{H}:\; \text{reversible lift}\;\longmapsto\CCfam{4}\quad \mathrm{D2-swap\_1}:\; \text{degree-two residue index}\;\longmapsto\CCfam{5}\\
&\mathrm{D2-swap\_2}:\; \text{degree-two residue index}\;\longmapsto\CCfam{6}\\
&\mathrm{D2-swap\_3}:\; \text{degree-two residue index}\;\longmapsto\CCfam{7}\\
&\mathrm{D2-fixed\_1}:\; \text{degree-two residue index}\;\longmapsto\CCfam{8}\\
&\mathrm{D2-fixed\_3}:\; \text{degree-two residue index}\;\longmapsto\CCfam{9}\\
&\mathrm{D2-fixed\_4}:\; \text{degree-two residue index}\;\longmapsto\CCfam{10}\\
&\mathrm{D2-fixed\_5}:\; \text{degree-two residue index}\;\longmapsto\CCfam{11}\\
&\mathrm{D2-fixed\_6}:\; \text{degree-two residue index}\;\longmapsto\CCfam{12}\\
&\mathrm{D2-fixed\_7}:\; \text{degree-two residue index}\;\longmapsto\CCfam{13}\\
&\mathrm{D3-I-01}\to\mathrm{D3-II-01}:\; k=1,e=1,m=0,n=2\;\longmapsto\CCfam{23}\\
&\mathrm{D3-I-02}:\; k=1,e=1,m=1,n=1\;\longmapsto\CCfam{14}\\
&\mathrm{D3-I-03}:\; k=1,e=1,m=2,n=0\;\longmapsto\CCfam{15}\\
&\mathrm{D3-I-04}:\; k=1,e=2,m=0,n=2\;\longmapsto\CCfam{16}\\
&\mathrm{D3-I-05}:\; k=1,e=2,m=1,n=1\;\longmapsto\CCfam{17}\\
&\mathrm{D3-I-06}:\; k=1,e=2,m=2,n=0\;\longmapsto\CCfam{18}\\
&\mathrm{D3-I-07}:\; k=2,e=1,m=0,n=1\;\longmapsto\CCfam{19}\\
&\mathrm{D3-I-08}:\; k=2,e=1,m=1,n=0\;\longmapsto\CCfam{20}\\
&\mathrm{D3-I-09}\to\mathrm{D3-II-01}:\; k=2,e=2,m=0,n=1\;\longmapsto\CCfam{23}\\
&\mathrm{D3-I-10}\to\mathrm{D3-II-02}:\; k=2,e=2,m=1,n=0\;\longmapsto\CCfam{24}\\
&\mathrm{D3-I-11}:\; k=3,e=1,m=0,n=0\;\longmapsto\CCfam{21}\\
&\mathrm{D3-I-12}:\; k=3,e=2,m=0,n=0\;\longmapsto\CCfam{22}\\
&\mathrm{D3-II-01}:\; m=0,n=1\;\longmapsto\CCfam{23}\quad \mathrm{D3-II-02}:\; m=1,n=0\;\longmapsto\CCfam{24}\\
&\mathrm{D3-III-01}:\; k=1,j=1,m=0,n=2\;\longmapsto\CCfam{25}\\
&\mathrm{D3-III-02}:\; k=1,j=1,m=1,n=1\;\longmapsto\CCfam{26}\\
&\mathrm{D3-III-03}:\; k=1,j=1,m=2,n=0\;\longmapsto\CCfam{27}\\
&\mathrm{D3-III-04}\to\mathrm{D3-III-03}:\; k=1,j=2,m=0,n=2\;\longmapsto\CCfam{27}\\
&\mathrm{D3-III-05}\to\mathrm{D3-III-02}:\; k=1,j=2,m=1,n=1\;\longmapsto\CCfam{26}\\
&\mathrm{D3-III-06}\to\mathrm{D3-III-01}:\; k=1,j=2,m=2,n=0\;\longmapsto\CCfam{25}\\
&\mathrm{D3-III-07}:\; k=2,j=1,m=0,n=1\;\longmapsto\CCfam{28}\\
&\mathrm{D3-III-08}:\; k=2,j=1,m=1,n=0\;\longmapsto\CCfam{29}\\
&\mathrm{D3-III-09}\to\mathrm{D3-III-08}:\; k=2,j=2,m=0,n=1\;\longmapsto\CCfam{29}\\
&\mathrm{D3-III-10}\to\mathrm{D3-III-07}:\; k=2,j=2,m=1,n=0\;\longmapsto\CCfam{28}\\
&\mathrm{D4-I-01}:\; k=1,e=1,m=0,n=1\;\longmapsto\CCfam{30}\\
&\mathrm{D4-I-02}:\; k=1,e=1,m=1,n=0\;\longmapsto\CCfam{31}\\
&\mathrm{D4-I-03}\to\mathrm{D2-fixed\_3}:\; k=1,e=2,m=0,n=1\;\longmapsto\CCfam{9}\\
&\mathrm{D4-I-04}\to\mathrm{D2-fixed\_4}:\; k=1,e=2,m=1,n=0\;\longmapsto\CCfam{10}\\
&\mathrm{D4-I-05}:\; k=1,e=3,m=0,n=1\;\longmapsto\CCfam{32}\\
&\mathrm{D4-I-06}:\; k=1,e=3,m=1,n=0\;\longmapsto\CCfam{33}\\
&\mathrm{D4-I-07}:\; k=2,e=1,m=0,n=0\;\longmapsto\CCfam{34}\\
&\mathrm{D4-I-08}\to\mathrm{D2-fixed\_5}:\; k=2,e=2,m=0,n=0\;\longmapsto\CCfam{11}\\
&\mathrm{D4-I-09}:\; k=2,e=3,m=0,n=0\;\longmapsto\CCfam{35}\\
&\mathrm{D4-II-01}:\; \text{degree-two residue index}\;\longmapsto\CCfam{36}\\
&\mathrm{D4-III-01}:\; j=1,m=0,n=1\;\longmapsto\CCfam{37}\\
&\mathrm{D4-III-02}:\; j=1,m=1,n=0\;\longmapsto\CCfam{38}\\
&\mathrm{D4-III-03}\to\mathrm{D2-swap\_1}:\; j=2,m=0,n=1\;\longmapsto\CCfam{5}\\
&\mathrm{D4-III-04}\to\mathrm{D2-swap\_1}:\; j=2,m=1,n=0\;\longmapsto\CCfam{5}\\
&\mathrm{D4-III-05}\to\mathrm{D4-III-02}:\; j=3,m=0,n=1\;\longmapsto\CCfam{38}\\
&\mathrm{D4-III-06}\to\mathrm{D4-III-01}:\; j=3,m=1,n=0\;\longmapsto\CCfam{37}\\
&\mathrm{D5-01}:\; e=1\;\longmapsto\CCfam{39}\quad \mathrm{D5-02}:\; e=2\;\longmapsto\CCfam{40}\quad \mathrm{D5-03}:\; e=3\;\longmapsto\CCfam{41}\\
&\mathrm{D5-04}:\; e=4\;\longmapsto\CCfam{42}\\
&\mathrm{D3-B-01}\to\mathrm{D3-I-02}:\; k=1,e=1,m=1,n=1\;\longmapsto\CCfam{14}\\
&\mathrm{D3-B-02}\to\mathrm{D3-I-03}:\; k=1,e=1,m=2,n=0\;\longmapsto\CCfam{15}\\
&\mathrm{D3-B-03}\to\mathrm{D3-I-05}:\; k=1,e=2,m=1,n=1\;\longmapsto\CCfam{17}\\
&\mathrm{D3-B-04}\to\mathrm{D3-I-06}:\; k=1,e=2,m=2,n=0\;\longmapsto\CCfam{18}\\
&\mathrm{D3-B-05}\to\mathrm{D3-I-08}:\; k=2,e=1,m=1,n=0\;\longmapsto\CCfam{20}\\
&\mathrm{D3-B-06}\to\!\mathrm{D3-B-09}\to\!\mathrm{D3-II-02}: k=2,e=2,m=1,n=0\;\longmapsto\!\CCfam{24}\\
&\mathrm{D4-B-01}\to\mathrm{D4-I-02}:\; k=1,e=1,m=1,n=0\;\longmapsto\CCfam{31}\\
&\mathrm{D4-B-02}\to\mathrm{D2-fixed\_4}:\; k=1,e=2,m=1,n=0\;\longmapsto\CCfam{10}\\
&\mathrm{D4-B-03}\to\mathrm{D4-I-06}:\; k=1,e=3,m=1,n=0\;\longmapsto\CCfam{33}\\
&\mathrm{D3-B-07}\to\mathrm{D3-III-02}:\; k=1,j=1,m=1,n=1\;\longmapsto\CCfam{26}\\
&\mathrm{D3-B-08}\to\!\mathrm{D3-B-07}\to\!\mathrm{D3-III-02}: k=1,j=2,m=1,n=1\;\longmapsto\!\CCfam{26}\\
&\mathrm{D3-B-09}\to\mathrm{D3-II-02}:\; m=1,n=0\;\longmapsto\CCfam{24}
\end{align*}
}

\paragraph*{The 27 original-parameter transport substitutions}
In each line the displayed source is contained in the target after the stated substitution. Variables on the right belong to the source; the variables assigned on the left belong to the target. Substitution makes every one of the ten original-coefficient differences identically zero. The source nonzero restrictions imply the target restrictions. The symbol $g$ remains the same original coefficient throughout.
\paragraph*{Transport 1: $\mathrm{D3-I-01}\to\mathrm{D3-II-01}$}
\noindent {\small $\displaystyle d=-35,\quad k=k$.}

\paragraph*{Transport 2: $\mathrm{D3-I-09}\to\mathrm{D3-II-01}$}
\noindent {\small $\displaystyle d=1,\quad k=k$.}

\paragraph*{Transport 3: $\mathrm{D3-I-10}\to\mathrm{D3-II-02}$}
\noindent {\small $\displaystyle d=0,\quad k=k$.}

\paragraph*{Transport 4: $\mathrm{D3-III-04}\to\mathrm{D3-III-03}$}
\noindent {\small $\displaystyle k=k$.}

\paragraph*{Transport 5: $\mathrm{D3-III-05}\to\mathrm{D3-III-02}$}
\noindent {\small $\displaystyle k=k,\quad u=u$.}

\paragraph*{Transport 6: $\mathrm{D3-III-06}\to\mathrm{D3-III-01}$}
\noindent {\small $\displaystyle k=k$.}

\paragraph*{Transport 7: $\mathrm{D3-III-09}\to\mathrm{D3-III-08}$}
\noindent {\small $\displaystyle k=k$.}

\paragraph*{Transport 8: $\mathrm{D3-III-10}\to\mathrm{D3-III-07}$}
\noindent {\small $\displaystyle k=k$.}

\paragraph*{Transport 9: $\mathrm{D4-I-03}\to\mathrm{D2-fixed\_3}$}
{\small
\begin{align*}
&a_{3}=\frac{640 g^{4} k}{27} + \frac{40 g^{4}}{27} - \frac{4 g^{2}}{3},\quad a_{4}=\frac{2560 g^{5} k}{81} + \frac{320 g^{5}}{243} - \frac{8 g^{3}}{27},\\
&a_{5}=- \frac{17920 g^{6} k}{243} + \frac{80 g^{6}}{243}.
\end{align*}
}

\paragraph*{Transport 10: $\mathrm{D4-I-04}\to\mathrm{D2-fixed\_4}$}
\noindent {\small $\displaystyle a_{3}=2 g^{4} k + \frac{13 g^{4}}{2} + 5 g^{2},\quad a_{4}=4 g^{5} k - 2 g^{5} - 3 g^{3},\quad a_{5}=- 2 g^{6} k + \frac{g^{6}}{4}$.}

\paragraph*{Transport 11: $\mathrm{D4-I-08}\to\mathrm{D2-fixed\_5}$}
\noindent {\small $\displaystyle a_{4}=- \frac{1344 g^{5}}{625} - \frac{56 g^{3}}{25},\quad a_{5}=- \frac{4368 g^{6}}{15625} - \frac{128 g^{4}}{625}$.}

\paragraph*{Transport 12: $\mathrm{D4-III-03}\to\mathrm{D2-swap\_1}$}
{\small
\begin{align*}
&a_{3}=- 384 g^{4} k - 1320 g^{4} + 84 g^{2},\quad a_{4}=- 1536 g^{5} k + 2880 g^{5} - 360 g^{3},\\
&a_{5}=10752 g^{6} k + 10800 g^{6}.
\end{align*}
}

\paragraph*{Transport 13: $\mathrm{D4-III-04}\to\mathrm{D2-swap\_1}$}
{\small
\begin{align*}
&a_{3}=- 384 g^{4} k - 1320 g^{4} + 84 g^{2},\quad a_{4}=- 1536 g^{5} k + 2880 g^{5} - 360 g^{3},\\
&a_{5}=10752 g^{6} k + 10800 g^{6}.
\end{align*}
}

\paragraph*{Transport 14: $\mathrm{D4-III-05}\to\mathrm{D4-III-02}$}
\noindent {\small $\displaystyle k=k$.}

\paragraph*{Transport 15: $\mathrm{D4-III-06}\to\mathrm{D4-III-01}$}
\noindent {\small $\displaystyle k=k$.}

\paragraph*{Transport 16: $\mathrm{D3-B-01}\to\mathrm{D3-I-02}$}
\noindent {\small $\displaystyle k=0,\quad u=v$.}

\paragraph*{Transport 17: $\mathrm{D3-B-02}\to\mathrm{D3-I-03}$}
\noindent {\small $\displaystyle k=0$.}

\paragraph*{Transport 18: $\mathrm{D3-B-03}\to\mathrm{D3-I-05}$}
\noindent {\small $\displaystyle k=0,\quad u=v$.}

\paragraph*{Transport 19: $\mathrm{D3-B-04}\to\mathrm{D3-I-06}$}
\noindent {\small $\displaystyle k=0$.}

\paragraph*{Transport 20: $\mathrm{D3-B-05}\to\mathrm{D3-I-08}$}
\noindent {\small $\displaystyle k=0$.}

\paragraph*{Transport 21: $\mathrm{D3-B-06}\to\mathrm{D3-II-02}$}
\noindent {\small $\displaystyle d=0,\quad k=0$.}

\paragraph*{Transport 22: $\mathrm{D4-B-01}\to\mathrm{D4-I-02}$}
\noindent {\small $\displaystyle k=0$.}

\paragraph*{Transport 23: $\mathrm{D4-B-02}\to\mathrm{D2-fixed\_4}$}
\noindent {\small $\displaystyle a_{3}=\frac{13 g^{4}}{2} + 5 g^{2},\quad a_{4}=- 2 g^{5} - 3 g^{3},\quad a_{5}=\frac{g^{6}}{4}$.}

\paragraph*{Transport 24: $\mathrm{D4-B-03}\to\mathrm{D4-I-06}$}
\noindent {\small $\displaystyle k=0$.}

\paragraph*{Transport 25: $\mathrm{D3-B-07}\to\mathrm{D3-III-02}$}
\noindent {\small $\displaystyle k=0,\quad u=v$.}

\paragraph*{Transport 26: $\mathrm{D3-B-08}\to\mathrm{D3-III-02}$}
\noindent {\small $\displaystyle k=0,\quad u=v$.}

\paragraph*{Transport 27: $\mathrm{D3-B-09}\to\mathrm{D3-II-02}$}
\noindent {\small $\displaystyle d=d,\quad k=0$.}

The eight duplicate removals, ten whole-family inclusions and nine boundary mergers are distinguished by the source assignments and transports above. The $k_0=0$ strata are included in the final maps. A type-II intersection with a quadratic family is a specialization, and does not remove the rest of that type-II family. This is why final coverage requires the displayed source transports, rather than pairwise intersections alone.

\endgroup

\end{document}